\documentclass{article}

\usepackage[T1]{fontenc}
\usepackage[english]{babel}
\usepackage{booktabs}
\usepackage{amsmath}
\usepackage{amsfonts}
\usepackage{amsthm}
\usepackage{amssymb}
\usepackage{color}
\usepackage{appendix}
\usepackage{subfigure}
\usepackage{array}
\usepackage{enumitem}
\usepackage{dcolumn}
\usepackage{siunitx}
\usepackage{verbatim}
\usepackage{xcolor}
\usepackage{authblk}

\usepackage[numbers]{natbib}
\usepackage{hypernat}
\usepackage[urlcolor=blue,colorlinks=true,citecolor=blue!50!black,linkcolor=black,plainpages=false,pdfpagelabels=true]{hyperref}%
\usepackage[all]{hypcap}

\usepackage{graphicx}
\graphicspath{{.}{./figures/}{./multiwave/figures/}}
\renewcommand{\leq}{\leqslant} 
\renewcommand{\geq}{\geqslant}

\renewcommand{\epsilon}{\varepsilon} 
\renewcommand{\hat}{\widehat}

\renewcommand{\tilde}{\widetilde}

\def\1{\mbox{1\hspace{-.35em}1}}
\def\R{\mathbb{R}}

\def\N{\mathbb{N}}

\def\E{\mathbb{E}}
\def\L{\mathbb{L}}

\def\Z{\mathbb{Z}}

\newcommand{\bff}{\mbox{${\boldsymbol{f}}$}}

\newcommand{\bd}{\mbox{${\mathbf d}$}}

\newcommand{\bI}{\mbox{${\mathbf I}$}}
\newcommand{\bu}{\mbox{${\mathbf u}$}}

\newcommand{\bA}{\mbox{${\mathbf A}$}}
\newcommand{\bB}{\mbox{${\mathbf B}$}}

\newcommand{\bG}{\mbox{${\mathbf G}$}}

\newcommand{\bX}{\mbox{${\mathbf X}$}}

\newcommand{\bW}{\mbox{${\mathbf W}$}}

\newcommand{\bOmega}{\mbox{\boldmath$\Omega$}}
\newcommand{\bLambda}{\mbox{\boldmath$\Lambda$}}

\newcommand{\bSigma}{\mbox{\boldmath$\Sigma$}}

\newcommand{\bepsilon}{\mbox{\boldmath$\epsilon$}}

\newcommand{\diag}{\mbox{diag}}

\newcolumntype{d}[1]{D{.}{.}{#1} } 

\newcommand{\argmin}{\displaystyle \mathop{argmin}}

\usepackage{listings}

\newcommand\pkg[1]{{\fontseries{b}\selectfont #1}}
\newcommand\code[1]{\texttt{#1}}
\newenvironment{CODE}{}{}
\lstnewenvironment{CodeIn}{\renewcommand\lstset{caption={}}}{}
\newenvironment{CodeOut}{\verbatim}{\endverbatim}

\providecommand{\keywords}[1]{\textbf{{Keywords.}} #1}

\author[a]{Sophie Achard}
\affil[a]{Univ. Grenoble Alpes, CNRS, Grenoble INP, GIPSA-lab, 38000 Grenoble, France
}
\author[b]{Ir{\`e}ne Gannaz}
\affil[b]{ Univ. Lyon, INSA de Lyon, CNRS UMR 5208, Institut Camille Jordan\protect\\
  20 avenue Albert Einstein, F-69621 Villeurbanne Cedex, France}

\title{Wavelet-based and Fourier-based multivariate Whittle estimation: {multiwave}}

\date{}

\begin{document}

\maketitle

\begin{abstract}
{Multivariate time series with long-dependence are observed in many applications such as finance, geophysics or neuroscience. Many packages provide estimation tools for univariate settings but few are addressing the problem of long-dependence estimation for multivariate settings. The package \pkg{multiwave} is providing efficient estimation procedures for multivariate time series. Two semi-parametric estimation methods of the long-memory exponents and long-run covariance matrix of time series are implemented. The first one is the Fourier-based estimation proposed by \cite{Shimotsu07} and the second one is a wavelet-based estimation described in \cite{AchardGannaz2014}. The objective of this paper is to provide an overview of the \textsf{R} package \pkg{multiwave} with its practical application perspectives.
}
\end{abstract}


%\vspace{\baselineskip}

\section{Introduction}

Time series with defined autocovariance functions are said to present long-memory or long-range dependency when their autocovariance function is decreasing very slowly, slower than an exponential decay. More precisely, let $g(\cdot)$ be the autocovariance function of a time series $X$. $X$ is said to be long-memory if there exists $\alpha$, $0<\alpha<1$, such that $g(t)$ is asymptotically equivalent to $|t|^{-\alpha}$ when $t \to +\infty$ (see \cite{beran.1994.1} and references therein). This definition implies that the covariance function is not summable. Equivalently, the spectral density $f(\cdot)$, if it exists, is such that $f(\lambda)$ is equivalent up to a constant to $|\lambda|^{1-\alpha}$ when $\lambda\to 0^+$. In this case, when trying to estimate the expectation using the empirical mean of long-memory time series, the variance of the estimator is not decreasing to 0 as $N^{-1}$ (where $N$ is the sample size). Hence it is crucial to take into account the presence of long-memory for defining good estimators \citep{beran.1994.1}. In the case of univariate time series, several very efficient approaches have been developed and validated. A web page entitled ``Time Series Analysis''\footnote{\url{https://cran.r-project.org/web/views/TimeSeries.html}} from the \textsf{R} software is providing a very exhaustive list of methods and softwares dealing with long-memory time series. Among others, we can cite the packages \pkg{fracdiff} \citep{haslett1989space}, \pkg{arfima} and \pkg{FGN} \citep{veenstra.2012.arfimaR}, \pkg{longmemo} \citep{beran.1994.1} and \pkg{forecast} \citep{hyndman.2008.JSS.forecastR}. For example, \pkg{fracdiff} is dedicated to simulation of fractional ARIMA time series and to estimation using regression of the periodogram. \pkg{longmemo} provides real data examples of time series with long-memory properties. 

 Approaches for multivariate long-memory time series are less developed. When dealing with multivariate time series, an important quantity to estimate is the covariance or correlation between pairs of time series. The effect of the presence of long-memory on this estimation is obvious, as stated by \cite{RobinsonDiscussion}. One \textsf{R} package, \pkg{waveslim}, is dedicated to the wavelet correlation analysis for pairs of random variables \citep{whitcher.2000.2} but long-range dependence properties are not considered. \cite{SelaHurvich2012} provide \textsf{R} code\footnote{freely available \url{http://pages.stern.nyu.edu/~rsela/VARFI/code.html}} for bivariate long-range dependent time series with parametric estimations. The objective of this paper is to provide an efficient \textsf{R} package, called \pkg{multiwave}, to estimate the long-memory parameters and covariance matrices for multivariate time series. The estimation procedures are based on a semi-parametric approach, which is robust to model misspecification. 
 
 The procedures are also suited for dealing with more than two dimensional data. Indeed they are based on Whittle approximation which provides a simple function to optimize. This function can be used for any dimension of the problem. In comparison, regression of the scalogram or periodogram \citep{Achard08} is based on a linear fit of pairs of time series, and thus there does not exist an easy way to extend to more than two dimensions.

 \pkg{multiwave} package is based on \cite{AchardGannaz2014}, where we developed a wavelet-based approach using Whittle approximation for an efficient estimation of the long-memory parameters and the long-run covariance matrices. In addition, \pkg{multiwave} proposes an implementation of an alternative method using Fourier decomposition as described in \cite{Shimotsu07}\footnote{code available in \textsf{matlab} \url{http://shimotsu.web.fc2.com/Site/Matlab_Codes.html}}.

As it is described in this paper, \pkg{multiwave} is a very versatile package and opens the way to estimation of the long-memory parameters and the long-run covariance matrices using multivariate data sets. It is in particular not necessary to assume the stationarity of the time series as it is the case when using Fourier decomposition \citep{FayMoulinesRoueffTaqqu}. The Whittle approximation is computed using either the coefficients of wavelet decomposition or the coefficients of Fourier decomposition when the time series are stationary.

The package \pkg{multiwave} is divided in three parts. A first group of functions is dedicated to the simulation of multivariate long-memory time series; the main function is \code{fivarma}. A second group of functions is implementing the wavelet decomposition, through \code{DWTexact} and associated functions. 
 Finally the computation of the estimators are coded using the Fourier decomposition in \code{mfw} and its derivatives 
  and using the wavelet decomposition in \code{mww} and its derivatives. 

The mathematical background is detailed in a separate Section \ref{sec:model}. The rest of the paper is dedicated to the description of the package \pkg{multiwave}. Simple examples of parametric models and real data are presented in Section \ref{sec:data} with corresponding functions of \pkg{multiwave} ready to apply. Core estimation functions using wavelets and Fourier transform are detailed in Section \ref{sec:estimation} along with pieces of code using simulated time series. Finally, practical considerations are discussed in the three last sections. Section \ref{sec:discussions} is discussing the practical choices of parameters. Comparaisons between wavelets and Fourier approaches are described in Section \ref{sec:simu}. And an application to real data in neuroscience is concluding the paper, Section \ref{sec:realdata}.    

\section{Theoretical background}
\label{sec:model}

As in the univariate case, the definition of long-memory for a $p$-vector process is based on the asymptotic behaviour of the cross-spectral density in the neighbourhood of zero \citep{Moulines07SpectralDensity}. We consider $N$ observations of a long-memory $p$-vector process $\bX=\{X_{\ell}(k),k\in\Z, \ell=1,\dots,p\}$, namely $\bX(1),\dots \bX(N)$. $\bX$ is said to be a multivariate $M(\bd)$ process when for each $\ell=1,\dots,p$ there exists $D_\ell\in\N$ such that the $D_\ell$-th order difference $\Delta^{D_\ell}X_\ell$ is covariance stationary. In addition, let us assume that for any $\ell,m=1,\dots,p$ the generalized cross-spectral density of $X_\ell$ and $X_m$ is
\begin{equation}
\label{eqn:density}
f_{\ell,m}(\lambda)=\frac{1}{2\pi}\Omega_{\ell,m}(1-e^{-i\lambda})^{-d_\ell}(1-e^{i\lambda})^{-d_m}f_{\ell,m}^S(\lambda),\qquad \lambda\in[-\pi,\pi],
\end{equation}
with $\bOmega=(\Omega_{\ell,m})_{\ell,m=1,\dots,p}$ an Hermitian matrix. The functions $f_{\ell,m}^S(\cdot)$ correspond to the short-memory dynamics of the process. The parameters $d_\ell$ satisfies $-1/2<d_\ell-D_\ell<1/2$.
More generally, the wavelet-based procedure is available for cross-spectral density satisfying an approximation
\begin{equation}\label{eqn:approx}\bff(\lambda)\sim  \bLambda(\bd) \bOmega \bLambda(\bd)^\ast,\quad \text{~~when~~}\lambda\to 0, \quad \text{~~with~~}  \bLambda(\bd)=\diag(|\lambda|^{-\bd}e^{-i \,\text{sign}(\lambda)\pi \bd/2}),
\end{equation}
where the exponent $\ast$ denotes the conjugate transpose operator. Here and subsequently $\sim$ means that the ratio of the left- and right-hand sides converges to one. Note that, the process $X_\ell$ is not necessarily stationary.

The long-range dependence parameter measures the power-like rate of decay of the autocovariance function. The long-run covariance matrix $\bOmega$ can be seen as the covariance at low frequencies between the time series. It gives a quantification of the link between the components of the multivariate time series. The long-run covariance parameter of the model is free from the difference in the autocorrelation behaviour of each component. It is linked with long-run correlations $\left(\Omega_{\ell,m}/\sqrt{\Omega_{\ell,\ell}\Omega_{m,m}}\right)_{\ell,m=1,\dots,p}$, which are also encountered in literature as power-law coherencies between two time series \citep{SelaHurvich2012} or as fractal connectivities \citep{Achard08}.

\subsection{A parametric example: FIVARMA}
\label{sec:fivarma}

Fractionally Integrated Vector Auto Regressive Moving Average (FIVARMA)   processes are parametric models with a spectral density satisfying approximation~\eqref{eqn:approx}. They correspond to Model A of \cite{Lobato97}. We refer to this paper for a detailed mathematical description.
 
Let $\bu$ be a $p$-dimensional white noise with $\E[\bu(t)\mid \mathcal F_{t-1}]=0$ and $\E[\bu(t)\bu(t)^T\mid \mathcal F_{t-1}]=\bSigma$, where $\mathcal F_{t-1}$ is the $\sigma$-field generated by $\{\bu(s),\,s<t\}$, and $\bSigma$ is a positive definite matrix.

Let $(\bA_k)_{k\in\N}$ be a sequence of $\R^{p\times p}$-valued matrices with $\bA_0$ the identity matrix and $\sum_{k=0}^\infty \|\bA_k\|^2<\infty$. The discrete Fourier transform of the sequence is denoted $\bA(\cdot)$, that is $\bA(\lambda)=\sum_{k=0}^{\infty} \bA_k e^{i k\lambda}$. We assume that all the roots of $|\bA(\L)|$ are outside the closed unit circle, where $\L$ denotes the lag operator. Let also $(\bB_k)_{k\in\N}$ be a sequence in $\R^{p\times p}$ with $\bB_0$ the identity matrix and $\sum_{k=0}^\infty \|\bB_k\|^2<\infty$. As defined for $\bA$, $\bB(\cdot)$ denotes the discrete Fourier transform of the sequence, $\bB(\lambda)=\sum_{k=0}^{\infty} B_k e^{i k\lambda}$.

Let $\bX$ be defined by

\begin{equation}
  \label{formula:cointegration}
  \bA(\L)\,\diag(\1-\L)^{\bd}\,\bX(t)=\bB(\L)\bu(t).
\end{equation}
  The spectral density satisfies

\[
f_{\ell,m}(\lambda)\sim_{\lambda\to 0^+} \frac{1}{2\pi}\Omega_{\ell,m}e^{-i\pi/2(d_\ell-d_m)}\lambda^{-(d_\ell+d_m)}
\] 
with \begin{equation}
\label{eqn:omega}
\bOmega=\bA(1)^{-1}\bB(1)\bSigma\bB(1)^T{\bA(1)^T}^{-1}.
\end{equation} Then $\bX$ is called a $FIVARMA(d,q)$ process and satisfies approximation~\eqref{eqn:approx}.

\subsubsection*{Limits of the model}

Note that in definition \eqref{formula:cointegration} the operators are applied in a given order, where the lag operator is taken first. Changing the order of the lag operator and autoregression corresponds to model $B$ of \cite{Lobato97} and $VARFI$ models of \cite{SelaHurvich2012} where $\bX$ is obtained with equation $\diag(\1-\L)^{\bd}\,\bA(\L)\,\bX(t)=\bB(\L)\bu(t)$, with similar notations than above.
That is, $\bX$ is obtained by fractional integration after autoregression, which is also called \textit{cointegration}. The spectral density still satisfies the approximation \eqref{eqn:approx} however the matrix $\bOmega$ may no longer be Hermitian.
Clearly \pkg{multiwave} package is not built to deal with such cases. We refer to alternative methods in literature, among others \cite{Robinson08, SelaHurvich2012, Shimotsu2012}. Taking into account cointegration is a difficult problem that exceeds the scope of this paper. Future work is needed to handle this particular case.

\subsection{Fourier-based estimation (MFW)}

The discrete Fourier transform and
the periodogram of $\bX$ evaluated at frequency $\lambda$ are defined as in \cite{Shimotsu07}'s procedure
\begin{align*}
\bW^F(\lambda)&= \frac{1}{\sqrt{2\pi N}} \sum_{t=1}^N \bX(t)e^{i t\lambda},\\
\bI^F(\lambda) &= \bW^F(\lambda)\bW^F(\lambda)^\ast.
\end{align*}

Let $\lambda_j = 2\pi j/N$, $j=1, \dots, m$, be the Fourier frequencies used in estimation, $m\in\N$. Define $\bLambda_j^F(\bd)=\diag\left(\lambda_j^{\bd}
e^{i(\pi-\lambda_j)\bd/2}\right)$.
The estimators $(\hat\bd^{MFW},\hat\bOmega^{MFW})$ are minimizers of the criterion
\begin{multline*}
\mathcal L^{MFW}(\bd,\bOmega)=\\ 
\frac{1}{m}\sum_{j=1}^{m} \left[ \log\det\left(\bLambda_j^F(\bd)\bOmega(\bd)
\bLambda_j^F(\bd)^\ast\right)+\bW^{F}(\lambda_j)^{\ast}\left(\bLambda_j^F(\bd)\bOmega(\bd)
\bLambda_j^F(\bd)^{\ast}
\right)^{-1}\bW^F(\lambda_j)\right].
\end{multline*}
The solution satisfies
\begin{align}
\label{eqn:dF}
\hat\bd^{MFW}&=\argmin_{\bd} \log\det(\hat \bOmega^{MFW}(\bd)) - 2\log(2)\left(\frac{1}{m}\sum_{j=1}^{m} \lambda_j\right),\\
\label{eqn:GF}
\hat \bOmega^{MFW}&=\hat \bOmega^{MFW}(\hat\bd^{MFW}),\\
\nonumber
\text{with~} \hat \bOmega^{MFW}(\bd)& =\frac{1}{m} \sum_{j=1}^{m} Re\left(\bLambda_j^F(\bd)^{-1}
\bI^F(j)\bLambda_j^F(\bd)^{-1}\right).
\end{align}

The dynamics of the frequencies at the neighbourhood of the origin is given by the dynamics of the spectral density around the zero frequency. The form of the criterion is justified by a second-order approximation of the spectral density matrix \eqref{eqn:density}, rather than the approximation \eqref{eqn:approx}. For an estimation based on the first-order approximation \eqref{eqn:approx}, one should replace $\bLambda_j^F(\bd)$ by $\bLambda_j^{F\,(1)}(\bd)=\diag\left(\lambda_j^{\bd}
e^{i\pi\bd/2}\right)$.

\cite{Shimotsu07} established the theoretical performance of this estimation procedure, for both the long-range dependence parameters and the long-run covariance matrix. It is shown that the variance for the estimation of the vector $\bd$ is decreased for the multivariate procedure with respect to a univariate one. It is worth mentioning that \cite{Lobato99} developed a similar estimation procedure, based on a rougher approximation of the cross-spectral density, $\bLambda_j^{F}(\bd)=\diag\left(\lambda_j^{\bd}\right)$. Interestingly, the quality of estimation for the vector $\bd$ is similar. Nevertheless, the estimation of the long-run covariance matrix $\bOmega$ is biased since it does not take into account the phase-shift appearing in $\bLambda_j^F(\bd)$. We refer to \cite{Lobato99} and to \cite{Shimotsu07} for a more detailed study of these estimators and their consistency.

\subsection{Wavelet-based estimation (MWW)}

Wavelets are providing a very efficient tool because of their high flexibility to deal with nonstationary time series which is particularly useful for real data applications. Their good performances in comparison to Fourier have already been shown for example in univariate settings  \citep{FayMoulinesRoueffTaqqu}. 

Let $(\phi(\cdot),\psi(\cdot))$ be respectively a father and a mother wavelets, satisfying regularity conditions, as stated in \cite{AchardGannaz2014}.

At a given resolution $j\geq 0$, for $k\in\Z$, we define the dilated and translated functions $\phi_{j,k}(\cdot)=2^{-j/2}\phi(2^{-j}\cdot -k)$ and $\psi_{j,k}(\cdot)=2^{-j/2}\psi(2^{-j}\cdot -k)$. The wavelet coefficients of the process $\bX$ are defined by $$\bW_{j,k}=\int_\R \tilde{\bX}(t)\psi_{j,k}(t)dt\quad j\geq 0, k\in\Z,$$ where $\tilde{\bX}(t)=\sum_{k\in\Z}\bX(k)\phi(t-k).$
For given $j\geq 0$ and $k\in\Z$, $\bW_{j,k}$ is a $p$-dimensional vector $\bW_{jk}=\begin{pmatrix}
W_{j,k}(1) & W_{j,k}(2) & \dots & W_{j,k}(p) \end{pmatrix}$ where $W_{j,k}(\ell)= \int_\R \tilde{\bX_\ell}(t)\psi_{j,k}(t)dt$.

For any $j\geq 0$, the process $(\bW_{j,k})_{k\in\Z}$ is covariance stationary \citep{AchardGannaz2014}. Let $\theta_{\ell,m}(j)$ denote the wavelet covariance at scale $j$ between processes $X_\ell$ and $X_m$, {\it i.e.} $\theta_{\ell,m}(j)=Cov(W_{j,k}(\ell), W_{j,k}(m))$ for any position $k$. Let us introduce the wavelet scalogram
\begin{equation}\label{eqn:scalogram}
\bI^W(j)=\sum_{k\in\Z} \bW_{j,k}\bW_{j,k}^T.
\end{equation}

The wavelet scalogram is the equivalent of the Fourier periodogram. Yet the scalogram is not normalized, on the contrary of the periodogram. We also introduce the function $K(\cdot)$, defined as 
\begin{equation}
\label{eqn:K}
K(\delta)=\int_{-\infty}^{\infty}{|\lambda|^{-\delta}|\hat\psi(\lambda)|^2\,d\lambda},\quad\delta\in(-\alpha,M).
\end{equation}

The wavelet Whittle procedure is described in \cite{AchardGannaz2014}. Let $\bLambda_j(\bd)=\diag\left(2^{j\bd},\, j_0\leq j\leq j_1\right)$. Let $\bG(\bd)$ denote a $p\times p$-matrix with $(\ell,m)$-{th} element equal to
\begin{equation}
\label{eqn:Gdef}
G_{\ell,m}(\bd)=f^S(0)\Omega_{\ell,m}K(d_\ell+d_m)cos({\pi(d_\ell-d_m)/2}).
\end{equation}
The estimators $(\hat\bd^{MWW},\hat\bG^{MWW})$ are defined by minimization of $\mathcal L^{MWM}(\bd,\bG)$, with 
\begin{multline*}
\mathcal L^{MWM}(\bd,\bG)=\\ \frac{1}{n}\sum_{j=j_0}^{j_1} \left[ n_j \log\det\left(\bLambda_j^W(\bd)\bG(\bd)\bLambda_j^W(\bd)\right)+\sum_k \bW_{j,k}^T\left(\bLambda_j^W(\bd)\bG(\bd)\bLambda_j^W(\bd)
\right)^{-1}\bW_{j,k}\right].
\end{multline*}
The estimation is here based on a first-order approximation of the spectral density matrix around 0.

 The solutions of this problem  satisfy
\begin{align}
\label{eqn:dhat}
\hat\bd^{MWW}&=\argmin_{\bd} \log\det(\hat \bG^{MWW}(\bd)) + 2\log(2)\left(\frac{1}{n}\sum_{j=j_0}^{j_1} j n_j\right)\left(\sum_{\ell=1}^p d_\ell\right),\\
\label{eqn:G}
\hat \bG^{MWW}(\bd) &=\frac{1}{n} \sum_{j=j_0}^{j_1} \bLambda_j^W(\bd)^{-1}
\bI^W(j)\bLambda_j^W(\bd)^{-1},
\end{align}
where $\bI^W(j)$ is the wavelet scalogram at scale $j$ defined in~\eqref{eqn:scalogram}. The long-run covariance matrix can then be estimated by
\begin{equation}
\label{eqn:Ohat}
\hat\Omega_{\ell,m}^{MWW}=\hat G_{\ell,m}^{MWW}(\hat \bd^{MWW})/(\cos(\pi(\hat d_\ell^{MWW}-\hat d_m^{MWW})/2) K(\hat d_\ell^{MWW}+\hat d_m^{MWW})).
\end{equation}

This second step in the estimation of the long-run covariance matrix $\bOmega$ is needed because the wavelets used in this paper are real and cannot correct the phase-shift (given by \eqref{eqn:Ohat}). This is not the case for the Fourier Whittle estimation described in \cite{Shimotsu07}.  Fortunately, the phase-shift can be expressed as a multiplicative cosine term in the covariance of the wavelet coefficients and a correction is still possible. 

\cite{AchardGannaz2014} established that the MWW estimators \eqref{eqn:dhat} and \eqref{eqn:Ohat} are consistent under non-restrictive conditions. The rate of convergence for the estimation of the long-range parameters $\bd$ is similar to the MFW estimator and is minimax. We refer to \cite{AchardGannaz2014} for the detailed study of the asymptotic behaviour of MWW estimation.

\section{Examples of multivariate long-memory time series}
\label{sec:data}

This section is describing specific functions of \pkg{multiwave} for the user to be able to simulate multivariate long-memory processes.
Parametric models are defined and implemented. In addition a data sets of real data from neuroimaging is provided. 

\subsection{Simulations of FIVARMA}

\label{sec:fivarma-code}

\pkg{multiwave} package proposes simulation functions for time series with a spectral density satisfying approximation \eqref{eqn:approx}. The main function is \code{fivarma} which computes a parametric FIVARMA process defined in Section~\ref{sec:fivarma}.

The input parameters of a $FIVARMA(q,d,r)$ process are the covariance matrix $\bSigma$ of the innovation process $\bu$, the vector AR (AutoRegressive) $(\bA_k)_{k=0,\dots,q}$, $\bA_k \in \R^{p\times p}$, the vector MA (MovingAverage) $(\bB_k)_{k=0,\dots,r}$, $\bB_k \in \R^{p\times p}$, and the vector of long-range parameters $\bd\in\R^p$.  FIVAR model of \cite{SelaHurvich2008} is a subcase, corresponding to MA coefficients equal to zero. The parameters of \code{fivarma} are thus, in order, \code{(N, d, cov\_matrix, VAR, VMA)} where \code{cov\_matrix}$=\bSigma$ and  \code{VAR} and \code{VMA} denote respectively the sequences of matrices  $(\bA_k)_{k=0,\dots,q}$ and $(\bB_k)_{k=0,\dots,r}$.

The output of the function \code{fivarma} is a list with first the values $\bX(1), \bX(2), \ldots, \bX(N)$ obtained by equation~\eqref{formula:cointegration} with $\bu(t)$ white noise with centered Gaussian distribution and covariance $\bSigma$. The second element of the list is the value of the matrix $\bOmega$ defined in \eqref{eqn:omega}.

\code{fivarma} is based on two other functions:
\begin{itemize}
\item \code{fracdiff} applies a vectorial fractional differencing procedure and corresponds to a FIVARMA(0,d,0).
\item \code{varma} computes a realisation of a multivariate ARMA process and corresponds to the case $\bd=0$.
\end{itemize}
Similar functions can be found in other packages ({\it e.g.} \pkg{fracdiff} and \pkg{MST} \citep{tsay2013multivariate}) but were re-implemented in \pkg{multiwave} package.

\vspace{\baselineskip}

\textit{Example.}
\begin{CODE}
\begin{CodeIn}
R> N <- 2^8
R> d0 <- c(0.2,0.4)
R> rho <- 0.8
R> cov <- matrix(c(1,rho,rho,1),2,2)
R> VMA <- diag(c(0.4,0.7))
R> VAR <- array(c(0.8,0.2,0,0.6),dim=c(2,2))
R> resp <- fivarma(N, d0, cov_matrix=cov, VAR=VAR, VMA=VMA)
R> x <- resp$x
R> long_run_cov <- resp$long_run_cov
R> long_run_cov
\end{CodeIn}
\begin{CodeOut}
          [,1]      [,2]
[1,] 0.6049383 0.5854938
[2,] 0.5854938 0.9730806
\end{CodeOut}
\end{CODE}
\begin{CodeIn}
R> par(mfrow=c(2,1),mai=c(0.5,1,0.5,0.5))
R> plot(x[,1],type='l',lty=1)
R> plot(x[,2],type='l',lty=1)
\end{CodeIn}

\includegraphics[width=10cm]{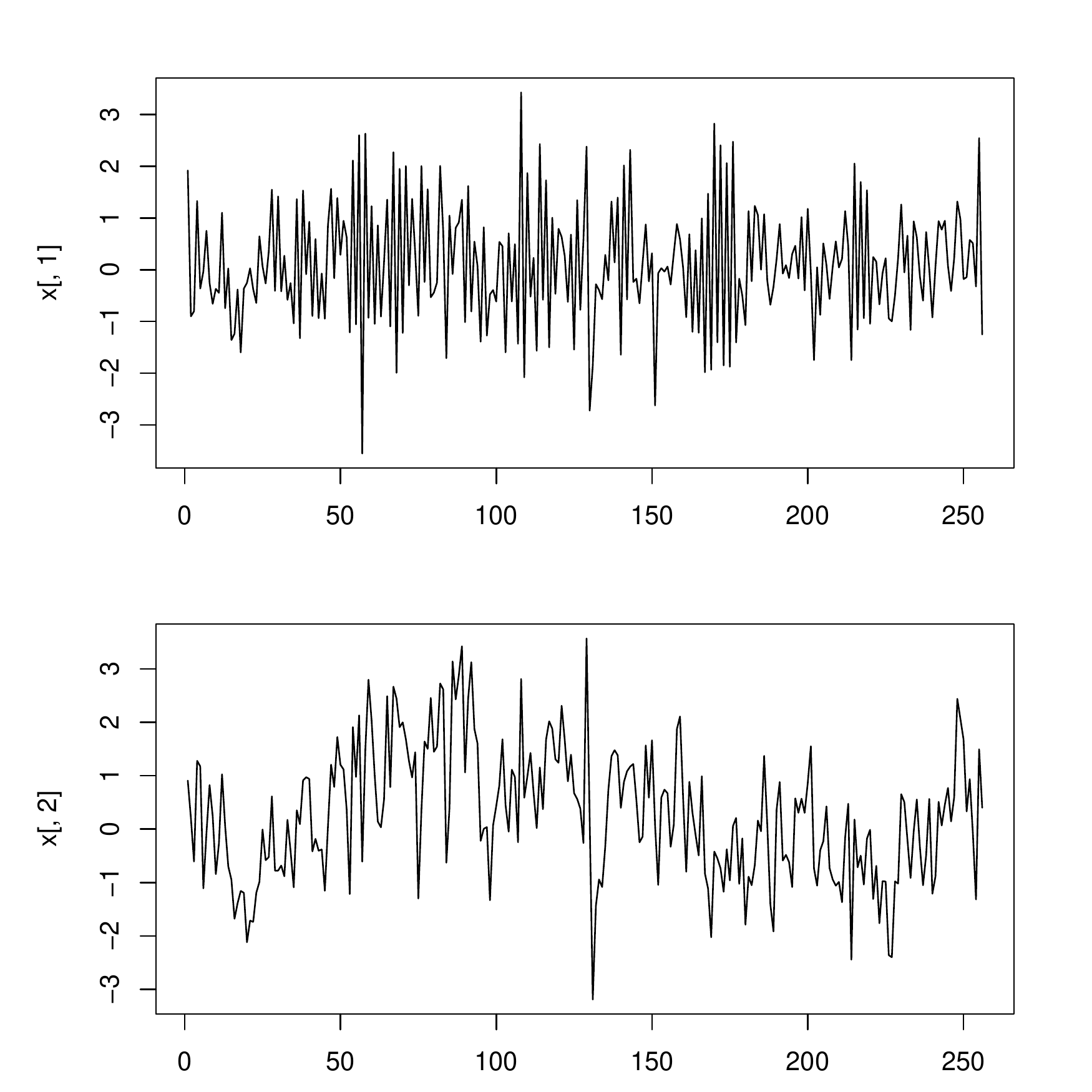}

\subsection{A real data set}
\label{sec:real_data_desc}
In order to describe how parameters can be chosen in a practical point of view, we provide a real data example (see Section~\ref{sec:realdata}).

\begin{figure}[!ht]
\caption{Plot of 6 arbitrary signals from a subject of fMRI data set.}
\label{fig:fmri.signal}
\begin{center}
{\includegraphics[scale=0.4]{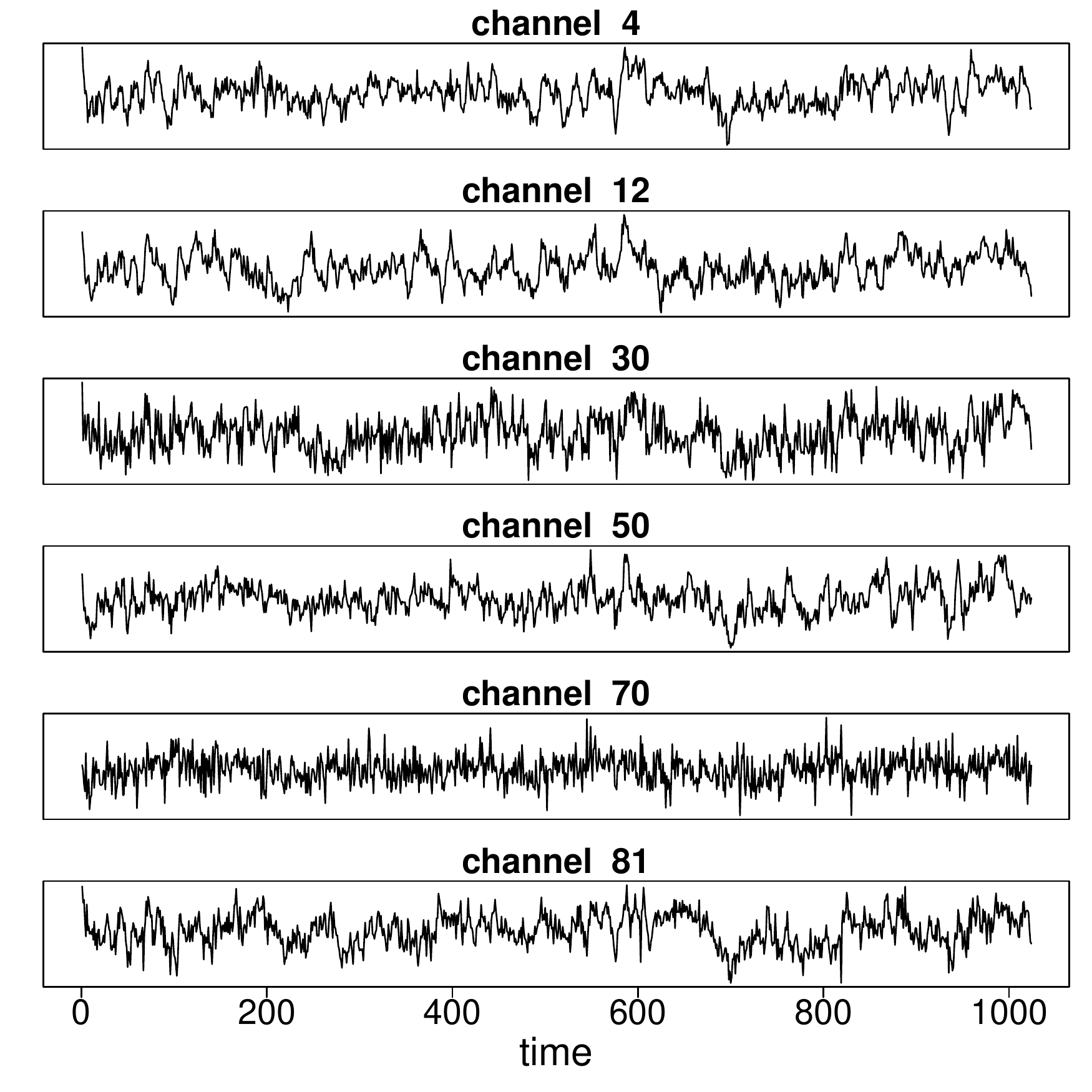}}
\end{center}
\end{figure}

Noninvasive data recorded from the brain are an example where the proposed methodology is efficient. The data consist of time series recording signals from the brain: electroencephalography (EEG) for the electrical signals, magnetoencephalography (MEG) for the magnetic signals or functional Magnetic Renonance Imaging (fMRI) for the Blood Oxygen Level Dependent (BOLD) signals. These data are intrinscally correlated because of the known interactions of the brain areas (also called regions of interest). Furthermore, it has already been shown that these time series present long-memory features \citep{maxim.2005.1}. Other data sets presenting similar features are coming from finance e.g. \cite{Songsiri}, where time series are correlated because of links between companies for example, and they also present long-memory characteristics. In this section, we observed time series extracted using fMRI facilities. The whole description of this data sets is detailed in \cite{Termenon2016}. 
The data set called \code{brainHCP} contains the time series of 1200 points in time and 89 regions of the brain. 
Figure \ref{fig:fmri.signal} displays 6 arbitrary signals from one subject in this data set.

\begin{CodeIn}
R> data(brainHCP)
R> dim(brainHCP)
\end{CodeIn}
\begin{CodeOut}
[1] 1200 89
\end{CodeOut}

\section{Estimation of the long-memory parameters and covariance matrix}
\label{sec:estimation}

The objective of this package is to provide the implementation of two sets of methods based either on Fourier or wavelet decomposition. That is  Multivariate Fourier Whittle (MFW) and Multivariate Wavelet Whittle (MWW) estimation procedures. The corresponding functions are respectively called \code{mfw} and \code{mww} in the package.

The output of the implemented methods consists in the estimation of two quantities, $\bd$ and $\bOmega$, where $\bd$ corresponds to the long-memory parameters of the time series and $\bOmega$ is reflecting the coupling between the pairs of time series.

\begin{figure}[!h]
\caption{Input and Output of \pkg{multiwave} package in a two-dimensional case.}
\begin{tabular}{m{7.5cm}m{2.2cm}m{0.4\textwidth}}
\begin{center}INPUT\end{center} & & \quad OUTPUT\\ 
\vspace{-0.5cm}\includegraphics[width=8cm]{plotx}&\begin{tabular}{c}
\pkg{multiwave}\\\includegraphics[width=2cm]{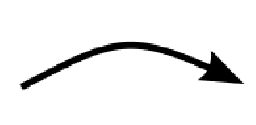}\end{tabular}
&\vline height 3cm depth 0.5cm \hspace{12pt}\begin{minipage}[b]{0.4\textwidth}
{Long-memory parameters

$\hat\bd = \begin{pmatrix}
0.17462\\ 0.24477
\end{pmatrix}$ 

\vspace{0.5cm}

Long-run covariance

$\hat\bOmega  = \begin{pmatrix}
 0.61064 & 0.55949 \\
 0.55949 & 1.29167
\end{pmatrix}
$
}
\end{minipage}
\end{tabular}

\end{figure}

The two functions \code{mfw} and \code{mww} are implementing the semiparametric Whittle estimation using respectively Fourier decomposition and wavelet decomposition in order to estimate $\bd$ and $\bOmega$.

A fast execution of these two functions is

\begin{CODE}
\begin{CodeIn}
R> ## Fourier decomposition
R> m <- 57 ## default value of Shimotsu 2017
R> res_mfw <- mfw(x,m)

R> ## Wavelet decomposition
R> res_filter <- scaling_filter('Daubechies',8) ## choice of filter
R> filter <- res_filter$h
R> LU <- c(2,11) ## choice of wavelet scales
R> res_mww <- mww(x,filter,LU)
\end{CodeIn}
%\begin{CodeOut}
%output of first input
%\end{CodeOut}
\end{CODE}
%$

\subsection{Multivariate Fourier Whittle estimation}
\label{sec:mfw}

As a first implementation, it is natural to use Fourier decomposition to approximate the spectral density of time series.

Package \pkg{multiwave} proposes functions to compute MFW estimators:
\begin{itemize}
\item \code{mfw} computes the multivariate Fourier Whittle estimators of both the long-range dependence parameters and the long-run covariance matrix.
\item \code{mfw\_cov\_eval} computes the multivariate Fourier-based Whittle estimator for the long-run covariance matrix for a given value of the long-range dependence $\bd$.
\item \code{mfw\_eval} returns the value of the multivariate Fourier Whittle criterion with respect to $\bd$ at a given value of $\bd$.
\end{itemize}

The functions \code{mfw\_cov\_eval} and \code{mfw\_eval} are internal functions of \code{mfw}. In \code{mfw}, we apply first a minimum search of \code{mfw\_eval} with respect to $\bd$, and \code{mfw\_cov\_eval} is returning the estimation of $\bOmega$ for the estimated value $\bd$.

We only detail function \code{mfw} hereafter and refer to the package description for other functions. 

Let  $\bX$ be the $p \times N$-matrix of observations, with general term ${x}_{\ell,i}=X_\ell(i)$, $\ell=1,\dots, p$ and $i=1,\dots, N$. Let $m$ be the number of frequencies used in MFW procedure. Given $x$ and $m$, the function \code{mfw} computes the MFW estimators defined by \eqref{eqn:dF} and \eqref{eqn:GF}, with the frequencies $\lambda_j = 2\pi j/N$, $j=1, \dots, m$. The optimization in equation \eqref{eqn:dF} is done using \code{optimize} function of \textsf{R} in one-dimensional settings and a Newton-type algorithm through \code{nlm} function of \textsf{R} otherwise. The initialization of the algorithm is set equal to the vector of univariate Fourier-based Whittle estimations. Even if it increases the computational time, such an initialization is important in high-dimensional settings. For example, in the MEG data set studied in \cite{AchardGannaz2014}, the optimization is done in $\R^{274}$. An initialization at the origin may not be able to reach the minimum even with a high number of iterations, due to the high dimension.

Function \code{mfw} returns a list with first the $p$-dimensional vector $\hat \bd^{MFW}$ and second the $p\times p$-matrix $\hat \bOmega^{MFW}$.  The quality of estimation is depending on the parameter $m$. Theoretical results show that $m$ must be small enough so the short-range properties of the time series do not bias the estimation. On the contrary a too small value will introduce variance in estimation since it decreases the number of frequency used in the procedure. \cite{Lobato99}, \cite{Shimotsu07} or \cite{Nielsen11} propose a default value $m=N^{0.65}$. This choice is discussed in the simulation study in Section \ref{sec:simu}.

\vspace{\baselineskip}

\textit{Example.}
\begin{CODE}
\begin{CodeIn}
R> # Simulation of the data
R> N <- 2^8
R> d0 <- c(0.2,0.4)
R> rho <- 0.8
R> cov <- matrix(c(1,rho,rho,1),2,2)
R> VMA <- diag(c(0.4,0.7))
R> VAR <- array(c(0.8,0.2,0,0.6),dim=c(2,2))
R> resp <- fivarma(N, d0, cov_matrix=cov, VAR=VAR, VMA=VMA)
R> x <- resp$x
R>
R> # Estimation
R> m <- N^(0.65) ## default value of Shimotsu
R> res_mfw <- mfw(x,m)
R> res_mfw
\end{CodeIn}
\begin{CodeOut}
$d
[1] 0.1239810 0.3497609

$cov
          [,1]      [,2]
[1,] 0.5004685 0.5840763
[2,] 0.5840763 1.2547056
\end{CodeOut}
\end{CODE}
%$

\subsection{Multivariate Wavelet Whittle estimation}
\label{sec:mww}

The functions applying MWW estimation in package \pkg{multiwave} are the following:
\begin{itemize}
\item \code{mww} computes the multivariate wavelet Whittle estimators of the long-range dependence parameters and the long-run covariance matrix.
\item \code{mww\_cov\_eval} computes the multivariate wavelet-based Whittle estimator for the long-run covariance matrix for a given value of the long-range dependence $\bd$.
\item \code{mww\_eval} returns the value of the multivariate wavelet Whittle criterion with respect to $\bd$ at a given value of $\bd$.
\end{itemize}
\code{mww\_cov\_eval} and \code{mww\_eval} are internal functions of \code{mww}. In \code{mww}, we apply first a minimum search of \code{mww\_eval} with respect to $\bd$, and \code{mww\_cov\_eval} is returning the estimation of $\bOmega$ for the estimated value of $\bd$. MWW estimation is based on the wavelet transform of time series and \code{mww} needs the definition of a wavelet filter. The computation of a filter and of a wavelet transform are described below.

\subsubsection{Wavelet transform and scalogram}

The wavelet decomposition in package \pkg{multiwave} is implemented using an exact discrete wavelet transform.
\begin{itemize}
\item  \code{scalingfilter} defines the wavelet filter (only Daubechies' wavelets are available).
\item  \code{computenj} computes the number of wavelet coefficients for each individual scale.
\item  \code{DWTexact} provides the wavelet transform of the data.
\item  \code{psi\_hat\_exact} gives the Fourier transform of the wavelet function.
\item \code{K\_eval} evaluates the value of the integral \eqref{eqn:K}.
\end{itemize}

\textit{Example.}
To obtain the wavelet filter of a Daubechies' wavelet of order 4, that is with 2 vanishing moments, one should write:
\begin{CODE}
\begin{CodeIn}
R> res_filter <- scaling_filter('Daubechies',4);
R> filter <- res_filter$h
R> filter
\end{CodeIn}
\begin{CodeOut}
[1]  0.4829629  0.8365163  0.2241439 -0.1294095
\end{CodeOut}
\end{CODE}
%$
Next, given an $N$ dimensional vector \code{x}, the wavelet coefficients of \code{x} are given by function \code{DWTexact}:
\begin{CODE}
\begin{CodeIn}
R> # Simulation of the data
R> N <- 2^8
R> d0 <- 0.2
R> resp <- fivarma(N, d0)
R> x <- resp$x
R>
R> # Wavelet decomposition
R> resw <- DWTexact(x,filter)
R> xwav <- resw$dwt # returns the vector of the wavelet coefficients
R> index <- resw$indmaxband 
R> Jmax <- resw$Jmax
R> Jmax
\end{CodeIn}
\begin{CodeOut}
[1] 6
\end{CodeOut}
\begin{CodeIn}
R> index
\end{CodeIn}
\begin{CodeOut}
[1] 127 189 219 233 239 241
\end{CodeOut}
\begin{CodeIn}
R> length(x)
\end{CodeIn}
\begin{CodeOut}
[1] 256
\end{CodeOut}
\end{CODE}
%$
\code{index} gives the index of the last coefficient at each scale and \code{Jmax} gives the maximal scale. The vector of coefficients \code{xwav} is $m$-dimensional, with $m$ maximum of \code{index}, equal to \code{index[Jmax]}. The coefficients of the third scale are for example given by:
\begin{CODE}
\begin{CodeIn}
R> xwav[seq(index[2]+1,index[3]),1]
\end{CodeIn}
\end{CODE}

Finally it is useful to compute the quantity $K(\delta)$ defined in equation \eqref{eqn:K}. Thus one needs to recover the Fourier transform of the wavelet, $\widehat{\psi}(\cdot)$. This is done using the function \code{psi\_hat\_exact}. Its inputs are the filter defined previously and an index of precision. It returns $(\widehat{\psi}(u_i))_{i=1,\dots,q*2^J}$ where $q$ is the length of the filter and $u_i$ are equally spaced points on the interval $[-\pi 2^{J-3} (q-1)/2; \,\pi 2^{J-3} (q-1)/2]$.
\begin{CODE}
\begin{CodeIn}
R> res_psi <- psi_hat_exact(filter,J=10)
R> psih <- res_psi$psih
R> gridh <- res_psi$grid
\end{CodeIn}
\end{CODE}
where \code{res\_psi\$grid} returns the values of the grid $(u_i)$ and \code{res\_psi\$psih} returns the corresponding values of $\widehat\psi(u_i)$. It is recommended to take $J\leq 15$ in practice and the default value is $J=10$. Indeed, a large value of \code{J} is increasing the computational time.\\
Given the function $\widehat{\psi}(\cdot)$, we are now able to evaluate $K(\bd)$ for a given value of $\bd$:
\begin{CODE}
\begin{CodeIn}
R> K <- K_eval(psih,gridh,d0)
\end{CodeIn}
\end{CODE}

\subsubsection{Estimation}

\code{mww} is now described in detail.
Let $\bX$ be the $p\times N$-matrix of observations, with general entries ${x}_{\ell,i}=X_\ell(i)$, $\ell=1,\dots, p$ and $i=1,\dots, N$. Let $LU$ be the bivariate vector giving the lowest scale $j_0$ and the upper scale $j_1$ of the wavelet coefficients used in estimation. Given $\bX$ and $LU$, the function \code{mww} computes the MWW estimators defined by \eqref{eqn:dhat} and \eqref{eqn:Ohat}. As previously, the optimization in equation \eqref{eqn:dhat} is done using \code{optimize} function of \textsf{R} in one-dimensional settings and a Newton-type algorithm through \code{nlm} function of \textsf{R} otherwise. The initialization of the algorithm is set equal to the vector of univariate wavelet-based Whittle estimations. The reasons are identical to the ones given for the function \code{mfw}. 

Function \code{mww} returns a list with first the $p$-dimensional vector $\hat \bd^{MWW}$ and second the $p\times p$-matrix $\hat \bOmega^{MWW}$.  The quality of estimation is depending on the parameters $j_0$ and $j_1$ appearing in \eqref{eqn:dhat} and \eqref{eqn:G}. Default value for $j_0$ is set to $2$ and for $j_1$ to the highest integer lower than $\log_2(N)$. The critical value to choose for estimation is $j_0$, as it can be seen in the theoretical conditions for consistency \citep{AchardGannaz2014} and in simulations studies. Similarly to the choice of the parameter $m$ for MFW procedure, a compromise exists between choosing a small value of $j_0$, which would introduce a bias due to the short-range properties of the time series, and a high value that will reduce the number of frequencies and thus increase variance.

\vspace{\baselineskip}

\textit{Example.}
\begin{CODE}
\begin{CodeIn}
R> # Simulation of the data
R> N <- 2^8
R> d0 <- c(0.2,0.4)
R> rho <- 0.8
R> cov <- matrix(c(1,rho,rho,1),2,2)
R> VMA <- diag(c(0.4,0.7))
R> VAR <- array(c(0.8,0.2,0,0.6),dim=c(2,2))
R> resp <- fivarma(N, d0, cov_matrix=cov, VAR=VAR, VMA=VMA)
R> x <- resp$x
R>
R> # Parameter of estimation
R> res_filter <- scaling_filter('Daubechies',8);
R> filter <- res_filter$h
R> LU <- c(2,8)
R>
R> # Estimation
R> res_mww <- mww(x, filter, LU)
R> res_mww
\end{CodeIn}
\begin{CodeOut}
$d
[1] 0.1955683 0.4857156

$cov
          [,1]      [,2]
[1,] 0.5826133 0.5051330
[2,] 0.5051330 0.9584876
\end{CodeOut}
\end{CODE}
%$

If one wants to apply several times the estimation on the same data set, or modifying the parameters of estimation, it is useful to separate the wavelet transform and the estimation scheme.
Wavelet-based estimation can be evaluated directly on the wavelet transform of the data using the following functions:
\begin{itemize}
\item \code{mww\_wav} computes the multivariate wavelet Whittle estimators of the long-range dependence parameters and the lon-run covariance matrix, given the wavelet transform of the data.
\item \code{mww\_wav\_cov\_eval} computes the MWW estimator for the long-run covariance matrix for a given value of the long-range dependence $\bd$, given the wavelet transform of the data.
\item \code{mww\_wav\_eval} returns the value of the multivariate wavelet-based Whittle criterion with respect to $\bd$ at a given value of $\bd$, for a specific wavelet transform of the data.
\end{itemize}
We refer to the description of the functions in the package for more details. 

\vspace{\baselineskip}

\textit{Example.}
\begin{CODE}
\begin{CodeIn}
R> # Simulation of the data
R> N <- 2^8
R> d0 <- c(0.2,0.4)
R> rho <- 0.8
R> cov <- matrix(c(1,rho,rho,1),2,2)
R> VMA <- diag(c(0.4,0.7))
R> VAR <- array(c(0.8,0.2,0,0.6),dim=c(2,2))
R> resp <- fivarma(N, d0, cov_matrix=cov, VAR=VAR, VMA=VMA)
R> x <- resp$x
R>
R> N <- dim(x)[1]
R> k <- dim(x)[2]
R>
R> # Parameter of estimation
R> res_filter <- scaling_filter('Daubechies',8)
R> filter <- res_filter$h
R> LU <- c(2,8)
R>
R> ## Wavelet decomposition
R> xwav <- matrix(0,N,k)
R> for(j in 1:k){
R> 	xx <- x[,j]	   
R>	resw <- DWTexact(xx,filter)
R>	xwav_temp <- resw$dwt
R>	index <- resw$indmaxband
R>	Jmax <- resw$Jmax
R>        xwav[1:index[Jmax],j] <- xwav_temp
R> }
R>	## we free some memory
R>	new_xwav <- matrix(0,min(index[Jmax],N),k)
R>	if(index[Jmax]<N){
R>		new_xwav[(1:(index[Jmax])),] <- xwav[(1:(index[Jmax])),]
R>	}
R>	xwav <- new_xwav
R>	index <- c(0,index)
R>
R> ##### Compute the wavelet functions
R> res_psi <- psi_hat_exact(filter,Jmax)
R> psih <- res_psi$psih
R> grid <- res_psi$grid
R>
R> # Estimation
R> res_mww_wav <- mww_wav(xwav,index,psih,grid,LU)
R> res_mww_wav
\end{CodeIn}
\begin{CodeOut}
$d
[1] 0.1684136 0.3693829

$cov
          [,1]     [,2]
[1,] 0.5177418 0.426430
[2,] 0.4264300 0.834756
\end{CodeOut}
\end{CODE}
%$

\section{Practical choices of parameters for MWW estimation}
\label{sec:discussions}

MWW procedure implemented in \code{mww} depends on mainly two parameters: the choice of the wavelet bases \code{filter} and the choice of the wavelet scales \code{LU}. These parameters are not fixed in the package as the performances of the estimation may be improved by a careful choice. The possibility to choose the wavelet scales is particularly of interest when dealing with short-range dependence. We show that a simple graphical representation of the scalogram is able to guide the user in the choice of the wavelet scales.

\subsection{Choice of the wavelet bases}
\label{sec:choixbase}

Actually \pkg{multiwave} only proposes Daubechies' wavelets, which satisfy theoretical properties of \cite{AchardGannaz2014}. Other bases are possible but not implemented. The wavelet bases is imputed {\it via} the parameter \code{filter} in function \code{mww}, where \code{filter} is obtained by \code{filter <- scaling\_filter('Daubechies',2*M)\$h}. %$ 
The main parameter characterizing the Daubechies' bases is thus the  number of vanishing moments, $M$. MWW estimation presents the advantage to be available even if the time series are nonstationary or with polynomial trends, as soon as $\max\bd\leq M$. For example, for stationary time series, a parameter $M=1$ is sufficient (which is equivalent to considering Haar bases). For real data application, when nonstationarity or trends are suspected, a higher value of $M$ is necessary.

As discussed in \cite{FayMoulinesRoueffTaqqu}, when $M$ increases, the quality of estimation (slightly) decreases. Depending on the data, a compromise is then needed between choosing a large enough number of vanishing moments $M$ to handle nonstationarity in the data and the quality of estimation.

On the contrary, MFW estimators are only suited to stationary time series. Some extensions of Fourier-based estimation were proposed in univariate setting such as tapered Fourier (see {\it e.g.} \cite{FayMoulinesRoueffTaqqu} and references therein). For multivariate estimation \cite{Nielsen11} proposes an extension of \cite{Shimotsu07} based on the transform defined in \cite{Abadir}. However, this approach gives satisfactory results only for $d < 1.5$ and we decided not to implement it in the package for simplicity.

\subsection{Choice of wavelet scales}
\label{sec:choixLU}

The second parameter we need to tune is \code{LU}, which corresponds to the range of scales used in estimation. \code{LU} is a two-dimensional vector, that is \code{LU<-c(j0, j1)}, with \code{j0} the lowest scale and \code{j1} the upper scale. Parameters \code{j0} and \code{j1} are respectively $j_0$ and $j_1$ defined in equations \eqref{eqn:dhat} and \eqref{eqn:G} in the estimation procedure.

One advantage of wavelets is to be able to qualitatively evaluate the choice of wavelet scales to estimate the long-memory parameters and correlation by inspection of wavelet scalogram. As mentioned in \cite{AbryVeitch98, FayMoulinesRoueffTaqqu} for univariate settings, the first and last scales may have to be discarded from the analysis. The first scale may be affected by the presence of short-memory phenomena. In the example of FIVARMA model, this is driven by the AR and MA coefficients. For the last scales, the impact is different and it comes from the finite length of the time series. Indeed, as derived in \cite{whitcher.2000.2}, the variance of the estimator is increasing with the wavelet scales. 

\begin{figure}[!h]
\caption[Boxplots of the log-variance of the wavelet coefficients at different scales for multivariate long-memory processes]{Boxplots of the log2-variance of the wavelet coefficients at different scales for simulated bivariate $FIVARMA(1,(0.4,0.4),1)$; (a) for the first component and (b) for the second component. (c) Boxplots of the log2-absolute covariance for the same data. The red lines represent the theoretical linear prediction given by $j (d_\ell+d_m)+\log_2(G_{\ell,m}(\bd))$, with $\ell=m=1$ for subfigure~(a), $\ell=m=2$ for subfigure~(b) and $\ell=1$, $m=2$ for subfigure~(c). The horizontal axis corresponds to increasing scales, that is, decreasing frequencies. The indexes of the horizontal axis display the number of coefficients available. %Indexes $i=1, 2$ in subfigure titles correspond to the component of the bivariate time series.
Parameters in FIVARMA model were the following: the white noise is Gaussian with a covariance equal to $\bSigma=\begin{pmatrix} 1 & 0.8 \\ 0.8 & 1\end{pmatrix}$, the AR coefficient is set equal to $A=\begin{pmatrix} 0.8 & 0 \\ 0.2 & 0.6\end{pmatrix}$ and the MA coefficient is set equal to $B=\begin{pmatrix} 0.4 & 0 \\0.2 & 0.7 \end{pmatrix}$.
 Calculation was done on $N=512$ observations for $100$ replications.}
\label{fig:scal}
\begin{center}{\includegraphics[width=1\textwidth]{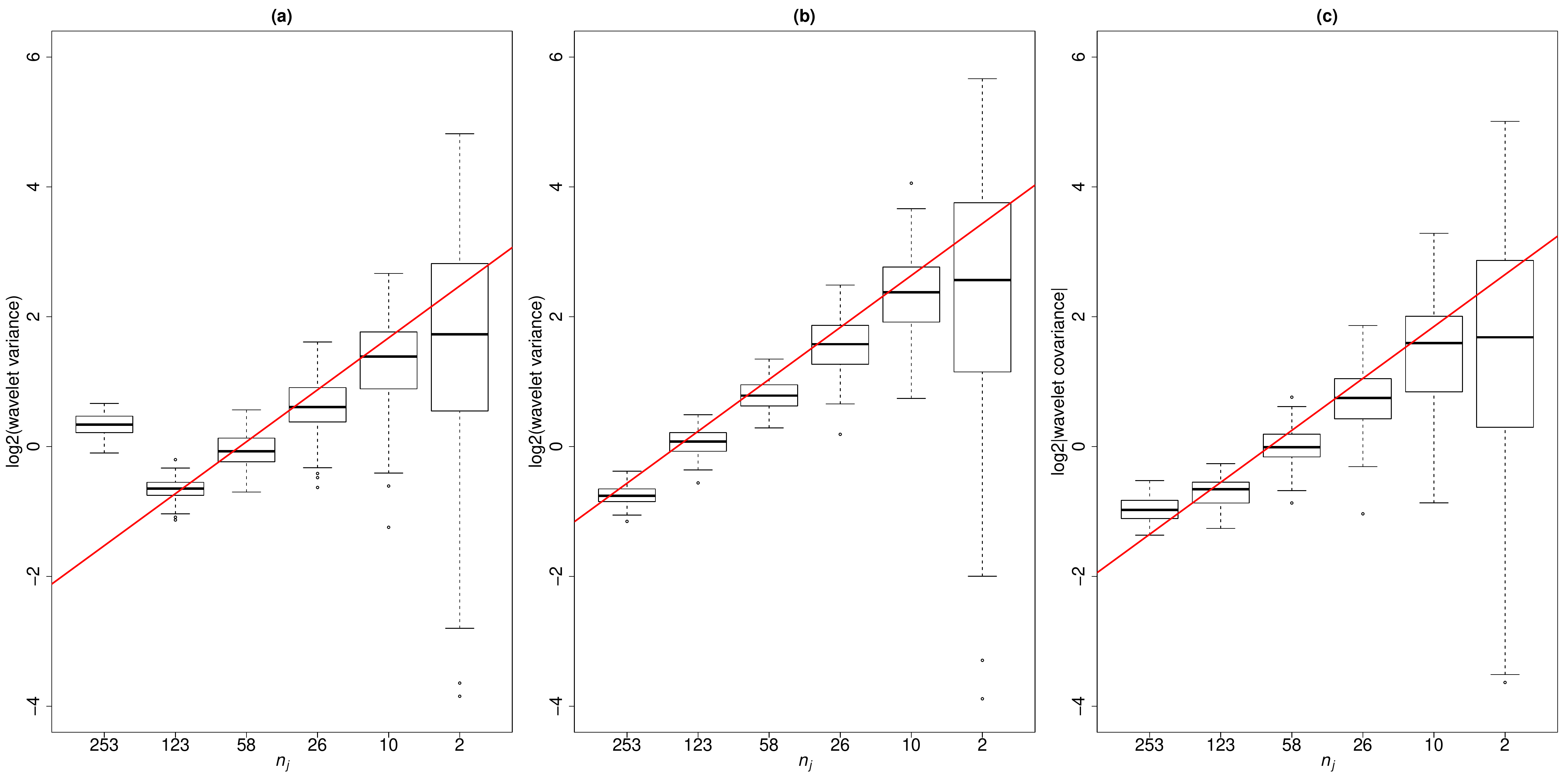}}
\end{center}
\end{figure}

The usual log-scalogram diagram used in univariate settings is showing the linear behaviour of the log variance with respect to the wavelet scales \citep{AbryVeitch98}. This is also true for the covariances as shown in Proposition~2 of \cite{AchardGannaz2014}. For all $k\in\Z$, $Cov(W_{j,k}(\ell),W_{j,k}(m))$ is equivalent to $2^{j(d_\ell+d_m)}G_{\ell,m}(\bd)$ when $j$ goes to infinity, with $G_{\ell,m}(\bd)$ defined in equation~\eqref{eqn:Gdef}. This property is illustrated in Figure~\ref{fig:scal} with a bivariate FIVARMA processes, as described in Section~\ref{sec:fivarma}. This figure represents the boxplots of the variance of the wavelet coefficients of each component of the time series at each scales (subfigures \ref{fig:scal}.(a) and \ref{fig:scal}.(b)). These plots correspond to the usual log-scalogram diagram \citep{AbryVeitch98}. Subfigure \ref{fig:scal}.(c) displays the analog representation of the covariance between the components. For both variance and covariance, the points satisfying the above approximation are aligned. Scales corresponding to non-aligned points should be removed from estimation as it can be seen on subfigure~\ref{fig:scal}.(a)~and~\ref{fig:scal}.(c) where the highest frequencies are modified by the presence of short-range dependence. Thus, one may discard the first scale from estimation to improve its quality.

\begin{figure}[!h]
\caption[Boxplots of the correlation of the wavelet coefficients at different scales for multivariate long-memory processes]{Boxplots of the correlation of the wavelet coefficients at different scales for four different simulated data: (a) bivariate $FIVARMA(0,(0.4,0.4),0)$ (b) bivariate $FIVARMA(1,(0.4,0.4),0)$; (c) bivariate $FIVARMA(1,(0.4,0.4),1)$; and finally, (d) an example with non stationary time series,  $FIVARMA(0,(0.8,1.2),0)$. The horizontal red lines represent the true long-run correlation for each simulation. The horizontal axis corresponds to increasing scales, that is, decreasing frequencies. The indexes of the horizontal axis display the number of coefficients available. 
Parameters in FIVARMA models were, if needed, the following: the white noise covariance is set equal to $\bSigma=\begin{pmatrix} 1 & 0.8 \\ 0.8 & 1\end{pmatrix}$, the AR coefficient is set equal to $A=\begin{pmatrix} 0.8 & 0 \\ 0.2 & 0.6\end{pmatrix}$ and the MA coefficient is set equal to $B=\begin{pmatrix} 0.4 & 0 \\0.2 & 0.7 \end{pmatrix}$.
Calculation was done on $N=512$ observations for $100$ replications.}
\label{fig:cov}
\begin{center}{\includegraphics[width=\textwidth]{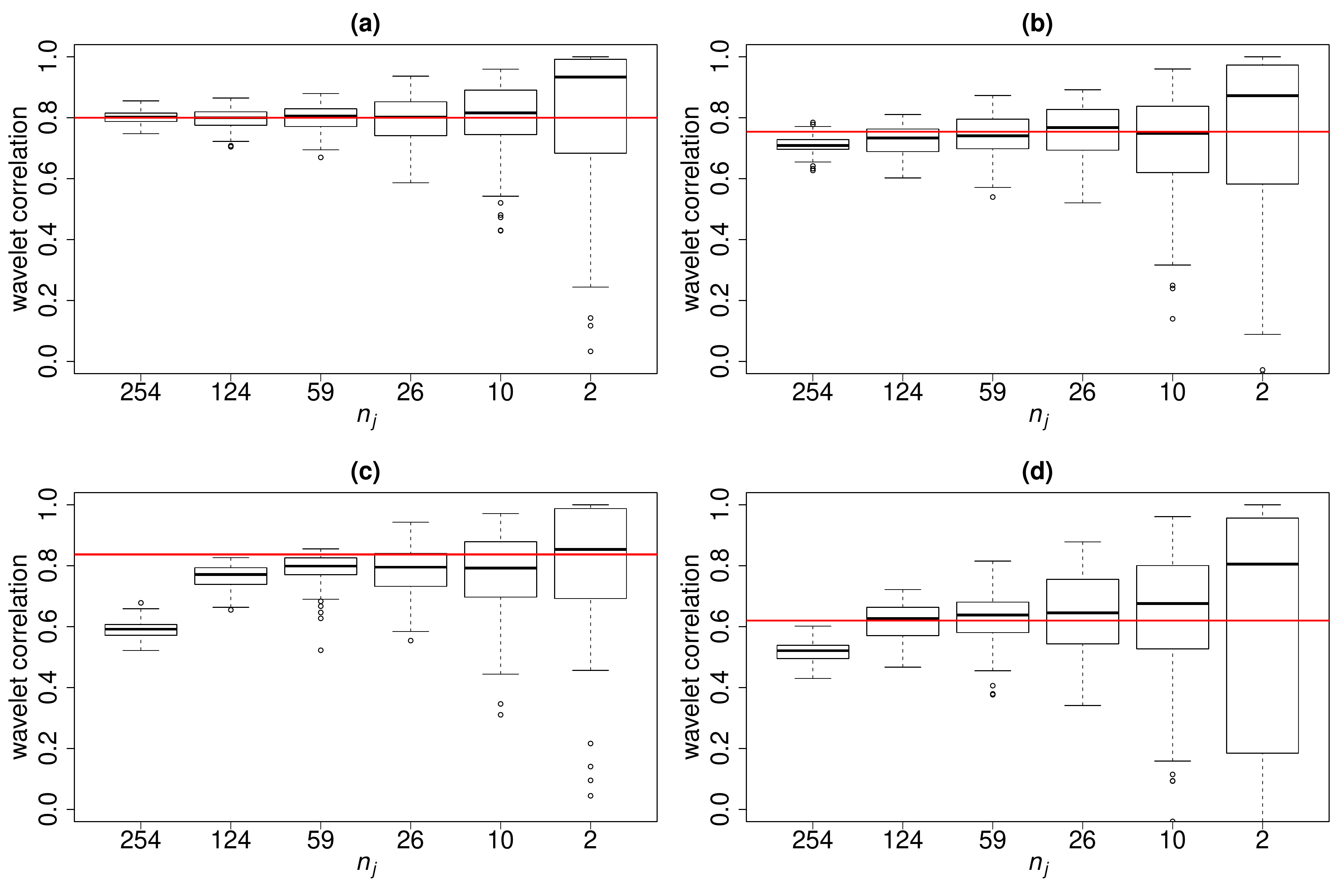}}
\end{center}
\end{figure}

Using the results on the wavelet variance and covariance in terms of scales, we showed that the wavelet correlation is asymptotically constant with respect to the wavelet scales. Indeed, for all $\ell,m=1,\dots,p$, for all $k\in\Z$, $Cor(W_{j,k}(\ell),W_{j,k}(m))$ is equivalent to $G_{\ell,m}(\bd)/\sqrt{G_{\ell,\ell}(\bd)G_{m,m}(\bd)}$ when $j$ goes to infinity \citep{AchardGannaz2014}. As for the log-scalogram diagram, the correlation between wavelet coefficients with respect to the scales can be plotted and scales where the observed correlation is not equal to the value obtained for the majority of scales should be removed from estimation. The wavelet correlation spectrum is a complementary way to qualitatively evaluate the range of scales where the analysis should be carried out. Figure \ref{fig:cov} illustrates this on four different data sets. Four different simulations of bivariate processes are applied using finite difference processes, FIVAR and FIVARMA processes, as described in Section~\ref{sec:fivarma}. Figure \ref{fig:cov} represents the boxplots of the correlation between wavelet coefficients of the two components of the time series at each scales. Again the presence of short-range dependence alters the highest frequencies (subfigures \ref{fig:cov}.(b) and \ref{fig:cov}.(c)). This is also observed with nonstationarity (subfigure \ref{fig:cov}.(d)). The first scales should then be removed from estimation. This free parameter of the package is particularly useful with the presence of short-range dependence or nonstationarity. Visual comparison to constant values may be easier for selection of the correct range of wavelet scale to use in the estimation. 

A similar discussion is detailed in Section~\ref{sec:realdata} for a real neuroscience data set. With real data sets a bootstrap procedure is necessary to obtain boxplots, as it will be explained in Section~\ref{sec:realdata}. For the Fourier procedure, the equivalent parameter is the number of frequencies $m$. However, wavelets are providing a graphical way to choose the upper and lower scales. To our knowledge, no equivalent qualitative evaluation for Fourier procedure is available.

\subsection{Numerical examples}

In order to quantify the quality of the choice of the parameters, numerical examples are provided. The first tables illustrate the quality of the estimators for a long-memory process with no short-range dependence. Then, the simulations are complexified by adding short-range behaviour or nonstationarity.

In each example we simulated 1000 Monte-Carlo replications of $N=512$ observations of FIVARMA(q,d,r), with a dimension $p=2$, for a set of different parameter values. Simulations are done using the function \code{fivarma}. MWW estimators are computed using function \code{mww}. MWW procedure is done using a Daubechies' wavelet bases with $M=4$ vanishing moments. This choice is motivated by discussion of Section~\ref{sec:choixbase}, because it can handle different settings, including nonstationary ones.

The quality of estimation is measured {\it via} the bias, the standard deviation ({std}) and the RMSE which is equal to $\sqrt{bias^2+std^2}$. For clarity, all tables are displayed at the end of the paper.

\subsubsection{A reference example}

We first consider a simple example, with neither short-range dependence, nor nonstationarity. Time series were simulated using a bivariate $FIVARMA(0,\bd,0)$ with a long-run correlation matrix $\bOmega=\begin{pmatrix}
1 & \rho \\ \rho & 1
\end{pmatrix}$ and $\rho =0.8$. The bivariate vector $\bd$ is chosen in $[0,0.5)^2$, such that the time series are stationary. 

As shown in Figure \ref{fig:cov}(a), and discussed in Section~\ref{sec:choixLU}, all scales can be kept for estimation.

Table~\ref{tab:dwav} displays results for the MWW estimation of $\bd$. This illustrates that multivariate estimation improves the quality of estimation for $\bd$. Indeed, the last column gives the ratio between the RMSEs of the multivariate wavelet-based estimation  and of the univariate wavelet-based estimation (ratio M/U). This ratio is always smaller than 1, that is, multivariate RMSE is always lower than univariate RMSE.

\subsubsection{Short-range dependence}

Consider a $FIVARMA(1,\bd,0)$ obtained with the model described in section \ref{sec:fivarma} and  given by the function \code{fivarma}. This case corresponds to a $FIVAR$ model of \cite{SelaHurvich2012}. The AR coefficient is taken equal to $\bA=\begin{pmatrix}
0.8 & 0 \\ 0.2 & 0.6
\end{pmatrix}$ and the correlation between the innovation processes equal to $\rho =0.8$. More precisely let $\bepsilon$ be a bivariate white noise process with covariance matrix $\bSigma=\begin{pmatrix}
1 & \rho \\ \rho & 1
\end{pmatrix}$ and let $\bu$ be the AR(1) process defined by $\bu(t)+\bA \bu(t-1)=\bepsilon(t)$. The time series observation $\bX(t)$ at time $t$ satisfies $(1-\L)^d\bX(t)= \bu(t).$ The matrix $\bOmega$ in equation \eqref{eqn:omega} is equal to
\[ \bOmega =(I+\bA)^{-1}\bSigma (I+\bA)^{-1\,T} \simeq \begin{pmatrix}
 0.3086 & 0.2392 \\
 0.2392 & 0.3260
\end{pmatrix}.
\]
The corresponding long-run correlation is thus equal to $0.754$.

As explained above, the finest scales are influenced by the short-range dependence and they have to be discarded from the estimation. As it can be seen in Figure~\ref{fig:cov}(c), the first two scales should be removed. We obtained accordingly that the lowest RMSE in estimation is obtained taking $j_0=3$. 

Numerical results for estimation are given in Table~\ref{tab:dwav2} and Table~\ref{tab:Owav2}. We can observe that estimation of $\hat\bd$ and $\hat\bOmega$ are still satisfactory. The RMSE is very similar to the previous case with no short-range behaviour.

\subsubsection{Nonstationarity}

Nonstationary examples are simulated using values of $\bd$ higher than $0.5$. We consider order 1 or 2 of nonstationarity, that is $\bd\in[0.5,2.5)^2$.

The behaviour of the wavelet correlations at each scale is illustrated in Figure \ref{fig:cov}(d). Contrary to stationary simulations where the optimal choice of $j_0$ was equal to $j_0=1$, the optimal choice of the parameter $j_0$ is $j_0=2$. Results are given in Table~\ref{tab:dwav3} and Table~\ref{tab:Owav3}.

Comparing Table~\ref{tab:dwav} and Table~\ref{tab:dwav3}, the quality of estimation of $\bd$ is still accurate in nonstationary settings, with similar values for the RMSE. As for the estimation of $\bOmega$, Table~\ref{tab:Owav} and Table~\ref{tab:Owav3} indicate that MWW still provides a good quality estimation of the long-run covariance matrix. The quality is slightly lower but still satisfactory.

\subsection{Discussion on identifiability}

\label{sec:identifiabilite}

In practical applications, it seems natural to assume that time series have the same order of stationarity. However, when two time series have long-memory parameters $d_\ell$ and $d_m$ satisfying $d_\ell-d_m=1$ the long-run covariance matrix $\bOmega$ is no longer identifiable with the wavelet-based procedure. Indeed, Proposition~2 in \cite{AchardGannaz2014} states that in this particular case the covariance $Cov(W_{j,k}(\ell),W_{j,k}(m))$ tends to 0 when the scale $j$ tends to infinity. Figure~\ref{fig:diff} illustrates this approximation for a bivariate $FIVARMA(0,(0.2,\,  1.2),0)$, a correlation matrix $\bOmega=\begin{pmatrix}
1 & \rho \\ \rho & 1
\end{pmatrix}$ and $\rho =0.8$.

When $\hat d_\ell- \hat d_m=1$, the estimator~\eqref{eqn:Ohat} is no longer defined. In practice, the quantity $\hat d_\ell- \hat d_m$ cannot be exactly equal to 1. Nevertheless, as dividing by a cosine function of this difference, a small error in the estimation of $(d_\ell, d_m)$ will lead to an important bias in the estimation of $\Omega_{\ell,m}$. As it can be seen in Figure~\ref{fig:bias.cov}, the resulting bias increases in the neighbourhood of the non-identifiable lines $d_\ell-d_m=\pm 1$.

{
\begin{figure}%[!h]
\caption[Boxplots of the covariance of the wavelet coefficients at different scales for a bivariate $FIVARMA(0,(0.2,1.2),0)$]{Boxplots of the covariance of the wavelet coefficients at different scales for a bivariate $ARFIMA(0,(0.2,1.2),0)$ with $\bOmega=\begin{pmatrix} 1 & 0.8 \\ 0.8 & 1\end{pmatrix}$. The index of the horizontal axis displays the number of coefficients available. Calculation was done on $N=512$ observations for $1000$ replications.}
\label{fig:diff}
\begin{center}{\includegraphics[width=12cm,height=7cm]{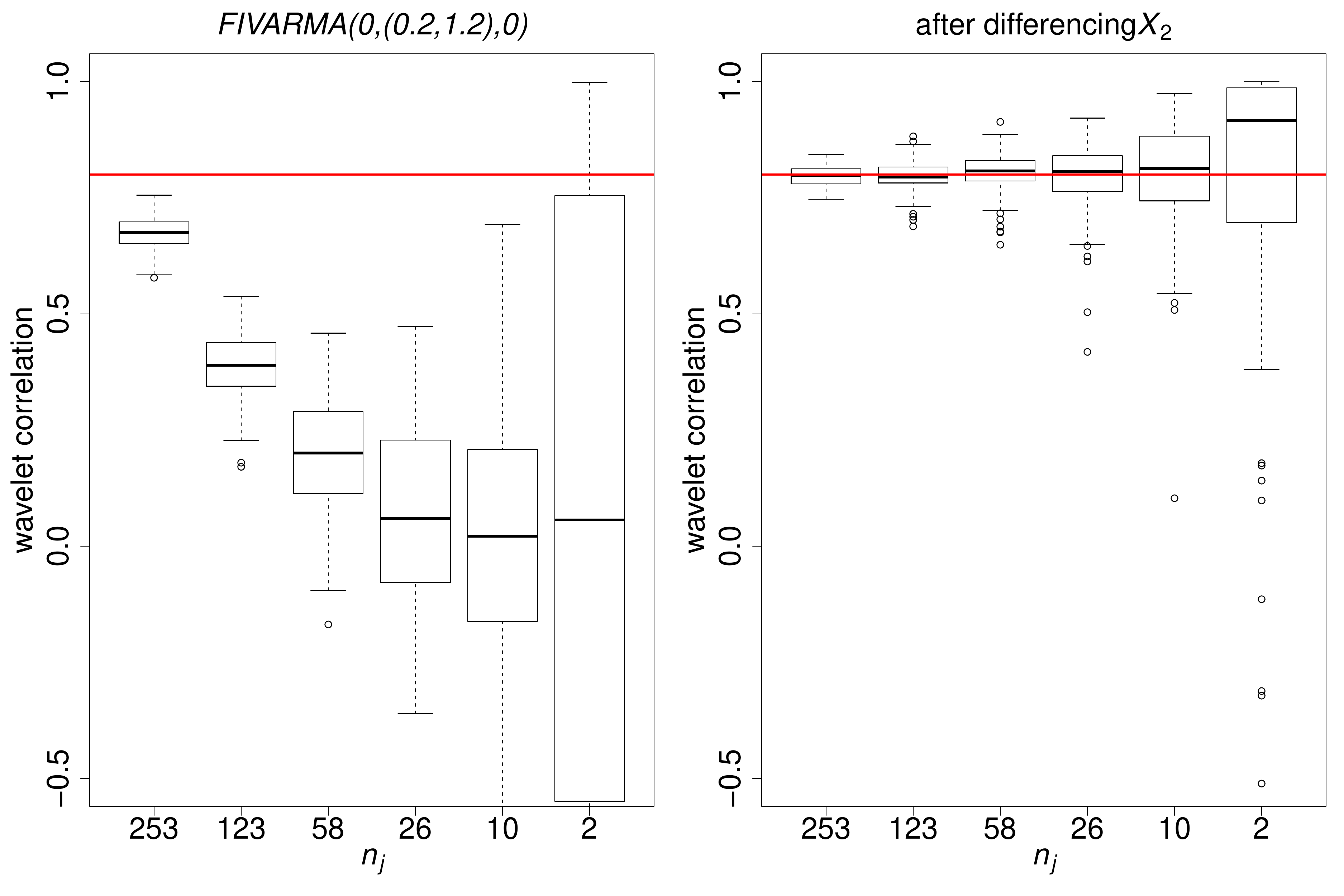}}
\end{center}
\end{figure}

\begin{figure}[!h]
\caption[RMSE in the estimation of the cross-covariance term $\Omega_{12}$ with respect to $(d_1,d_2)$.]{RMSE in the estimation of the cross-covariance term $\Omega_{12}$ with respect to $(d_1,d_2)$. Estimation was done using multivariate Wavelet Whittle estimator in a bivariate $FIVARMA(0,(d_1, d_2),0)$ with $\mathbf{\Omega}=\begin{pmatrix} 1 & 0.8\\ 0.8 & 1 \end{pmatrix}$. Subfigure~(b) represents an image plot of subfigure~(a) (with a different colour scale to improve visual quality). Blue lines on subfigure~(b) correspond to $d_2-d_1=\pm1$. Calculation was done on $N=512$ observations for $1000$ replications.}
\label{fig:bias.cov}
\begin{center}
\begin{tabular}{cc}
  (a) & (b)\\
 [-0cm]
  {\includegraphics[width=8cm]{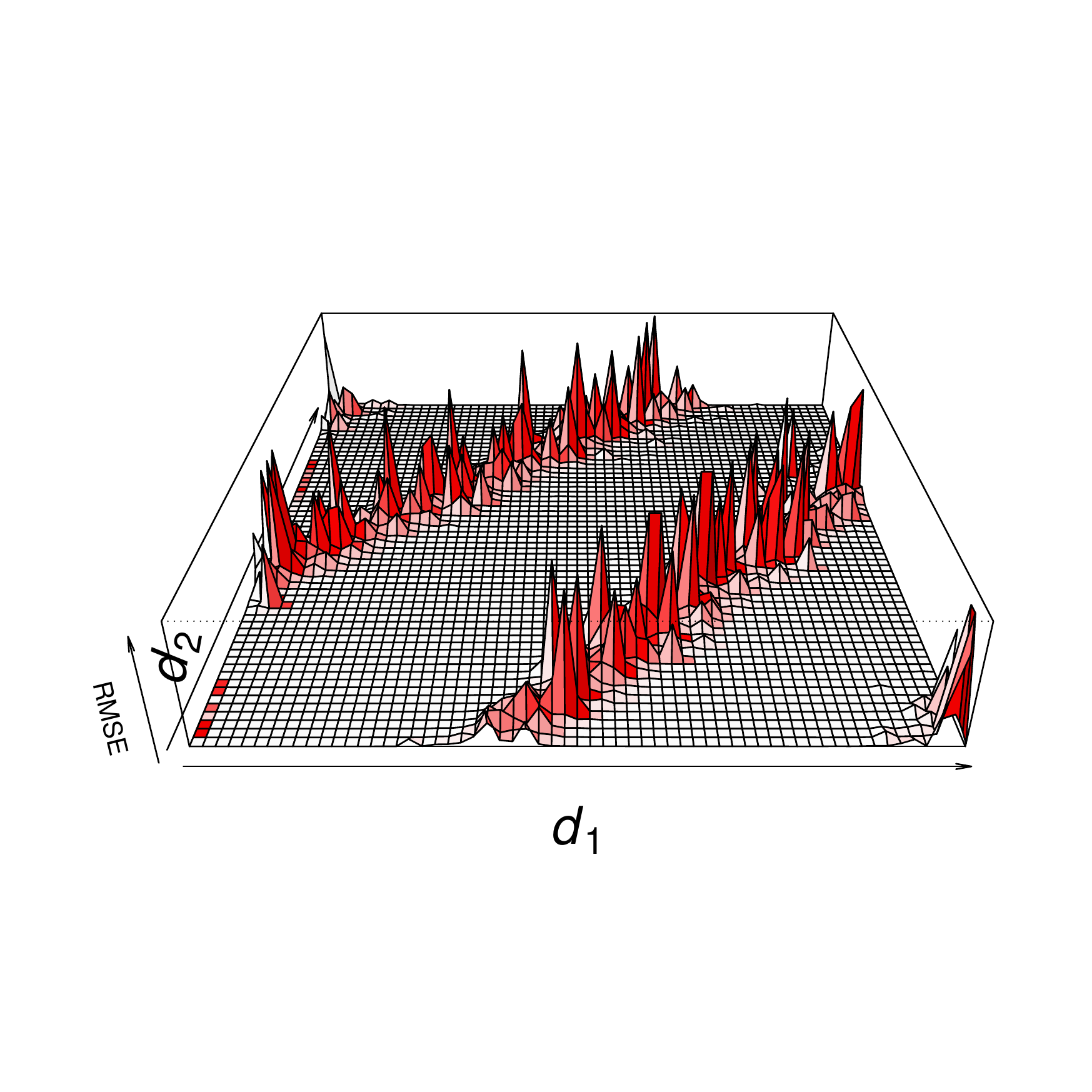}}&
\begin{minipage}[b]{0.4\textwidth}{\vspace*{1cm}\includegraphics[width=7cm]{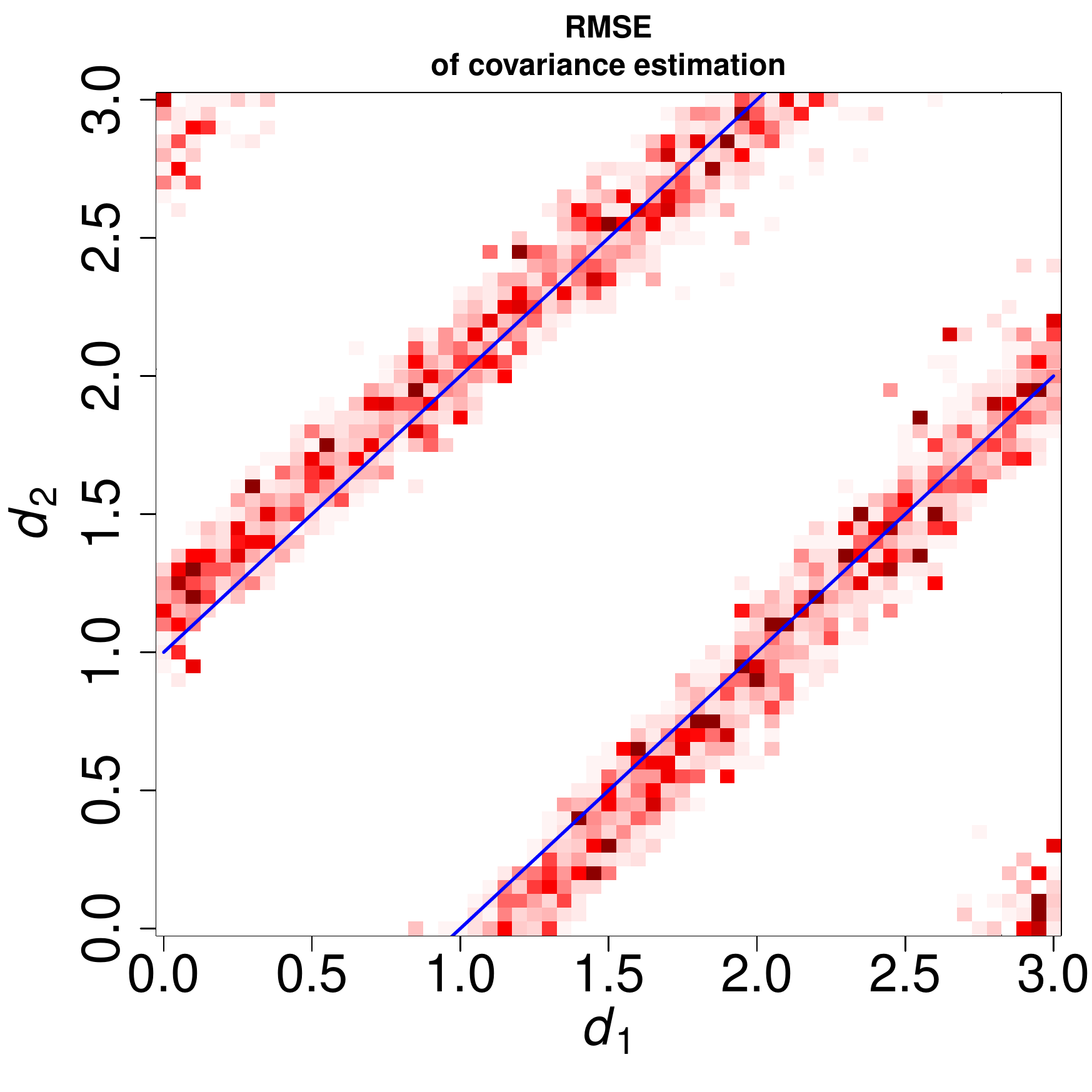}\vspace{1.5cm}}\end{minipage}
\vspace{-1cm}
\end{tabular}
\end{center}
\end{figure}
}

When this situation occurs, say when the difference between $d_\ell-d_m$ is between 0.75 and 1.25, the estimation of $\bd$ is not affected. But the user must be careful for the estimation of $\bOmega$. One solution is to differentiate or integrate one of the two processes. For example, Table~\ref{tab:bias.cov} illustrates the non-identifiability of $\bOmega$ in a bivariate $FIVARMA(0,\begin{pmatrix}
0.2 & 1.2 \end{pmatrix},0)$. When differentiating the second component (with $d_2=1.2$) the estimator has again good performances.

 \begin{table}
   \caption{Multivariate Wavelet Whittle estimation of $\bOmega$ for a bivariate $FIVARMA(0,(0.2,\,1.2),0)$ with $\rho=0.8$, $N=512$ with 1000 repetitions. As the memory parameter of the second component is greater than one, one possibility is to differentiate the second component. $j_0$ is chosen to be equal to 1. }
   \label{tab:bias.cov}
\begin{center}
\begin{tabular}{@{}lSSSSSSS@{}} %\hspace{0.5cm}
\toprule
& \multicolumn{3}{l}{\textbf{Without differentiation }} & & \multicolumn{3}{l}{\textbf{With differentiation} }\\
\cmidrule{2-4} \cmidrule{6-8}
$\bOmega$ & \multicolumn{1}{c}{\emph{bias}} & \multicolumn{1}{c}{\emph{std}} & \multicolumn{1}{c}{\emph{RMSE}} & & \multicolumn{1}{c}{\emph{bias}} & \multicolumn{1}{c}{\emph{std}} & \multicolumn{1}{c}{\emph{RMSE}} \\
\midrule 
 $\Omega_{1,1}$ &  0.0935 & 0.0762 & 0.1206 && 0.0349 & 0.0678 & 0.0762\\ 
 $\Omega_{1,2}$ &  4.1863 & 7.2103 & 8.3375 && 0.0261 & 0.06 & 0.0654\\ 
 $\Omega_{2,2}$ &  0.2215 & 0.0819 & 0.2362 && 0.0292 & 0.07 & 0.0758\\ 
 correlation   &  3.5255 & 6.2935 & 7.2137 &&  0.0003 & 0.0155 & 0.0155\\
\bottomrule
\end{tabular}
\end{center}
\end{table}

\section{MFW estimation and comparison with MWW}
\label{sec:simu}

The comparison between Fourier-based and wavelet-based approach is presented now. Time series were simulated using a bivariate $FIVARMA(0,\bd,0)$ with a long-run correlation matrix $\bOmega=\begin{pmatrix}
1 & \rho \\ \rho & 1
\end{pmatrix}$ and $\rho =0.8$. The bivariate vector $\bd$ is chosen in $[0,0.5)^2$, such that the time series are stationary. In such a setting, Fourier-based estimators are available.

  For Fourier-based approach, the parameter to choose is \code{m} corresponding to the number of frequencies taken into account in the estimation.  The default value in \cite{Shimotsu07} is $m=N^{0.65}$. We also make comparisons with an optimal value computed by minimizing the RMSE.

\subsection{Estimation of the long-memory parameters}

 Table~\ref{tab:dfou} gives the results obtained for the estimation of $\bd$ using MFW procedure (to be compared to Table~\ref{tab:dwav} for wavelets). Comparison between Fourier and wavelet-based procedures is summarized by the ratio between the RMSE given by MWW estimation and the RMSE given by MFW estimation, denoted by ratio W/F. Taking the same number of frequencies as \cite{Shimotsu07}, that is, $m=N^{0.65}$, Table~\ref{tab:dfou} shows that the quality of MWW end MFW procedures are comparable, even if wavelet-based estimation slightly improves Fourier-based estimation with such a choice of $m$.

Next we also consider the number of frequencies leading to the minimal RMSE for MFW estimation. As it can be seen in Table~\ref{tab:dfou}, qualities of both procedures are very similar but MFW estimation then (slightly) surpasses MWW estimation.

Very precise comparisons of Fourier-based and wavelet-based approaches are described in \cite{FayMoulinesRoueffTaqqu} for a univariate setting. In particular, it is shown that since the time series are stationary, the use of Haar bases should improve MWW quality. The authors indeed obtained better results with the Haar-based procedure than with Fourier-based procedure. To highlight the versatility of wavelet-based procedures, we choose here wavelet bases with four vanishing moments. A fair comparison with a Fourier-based method should consider tapered Fourier of order 4 as it is detailed in \cite{FayMoulinesRoueffTaqqu}. They quantify the influence of the regularity of the wavelet bases and discuss the comparison with (tapered) Fourier bases. Similar results are expected to be obtained in the case of multivariate time series, however, this topic exceeds the scope of this paper. 

\subsection{Estimation of the long-run covariance}

Finally, Table~\ref{tab:Ofou} and Table~\ref{tab:Ofou2} display results for the estimation of $\bOmega$ with MFW method (to be compared to Table~\ref{tab:Owav} for wavelets). When MFW estimation is applied with the usual number $m=N^{0.65}$ of frequencies, one can see that the wavelet-based procedure still estimates better the long-run covariance and the long-run correlation, with a ratio W/F always lower than $0.6$ for the estimation of $\bOmega$ terms. When the number of frequencies in Fourier-based estimation is chosen optimally, MFW and MWW procedures behave similarly and none appears significantly better than the other.

To conclude, MFW and MWW estimation procedures give very similar results. The slight improvement of Fourier-based procedure for the estimation of $\bd$ can be explained by the choice of the wavelet bases, however wavelets are efficient for a large set of applications, including time series with trends and nonstationarity features.

\section{Application on real neuroscience data}

\label{sec:realdata}

\begin{figure}[!ht]
\caption[Boxplots of the correlation of the wavelet coefficients at different scales for fMRI data]{Boxplots of the correlation of the wavelet coefficients at different scales for real time series from fMRI data sets: (a) Time series 1 and 2; (b) Time series 13 and 14; (c) Time series 31 and 32; (d) Time series 47 and 48. boxplots were obtained using sliding windows with $N=512$ points, extracted from two fMRI time series with length equal to 1200 points, from a single subject. The estimated long parameters $d$ of the two time series are equal. fMRI data set is described in section~\ref{sec:realdata}. The index of the horizontal axis displays the number of coefficients available. The horizontal red lines represent the estimated long-run correlation. Calculation was done on $N=512$ observations for $100$ replications using sliding windows (with overlap).}
\label{fig:cov_real_data}
\begin{center}{\includegraphics[width=0.92\textwidth]{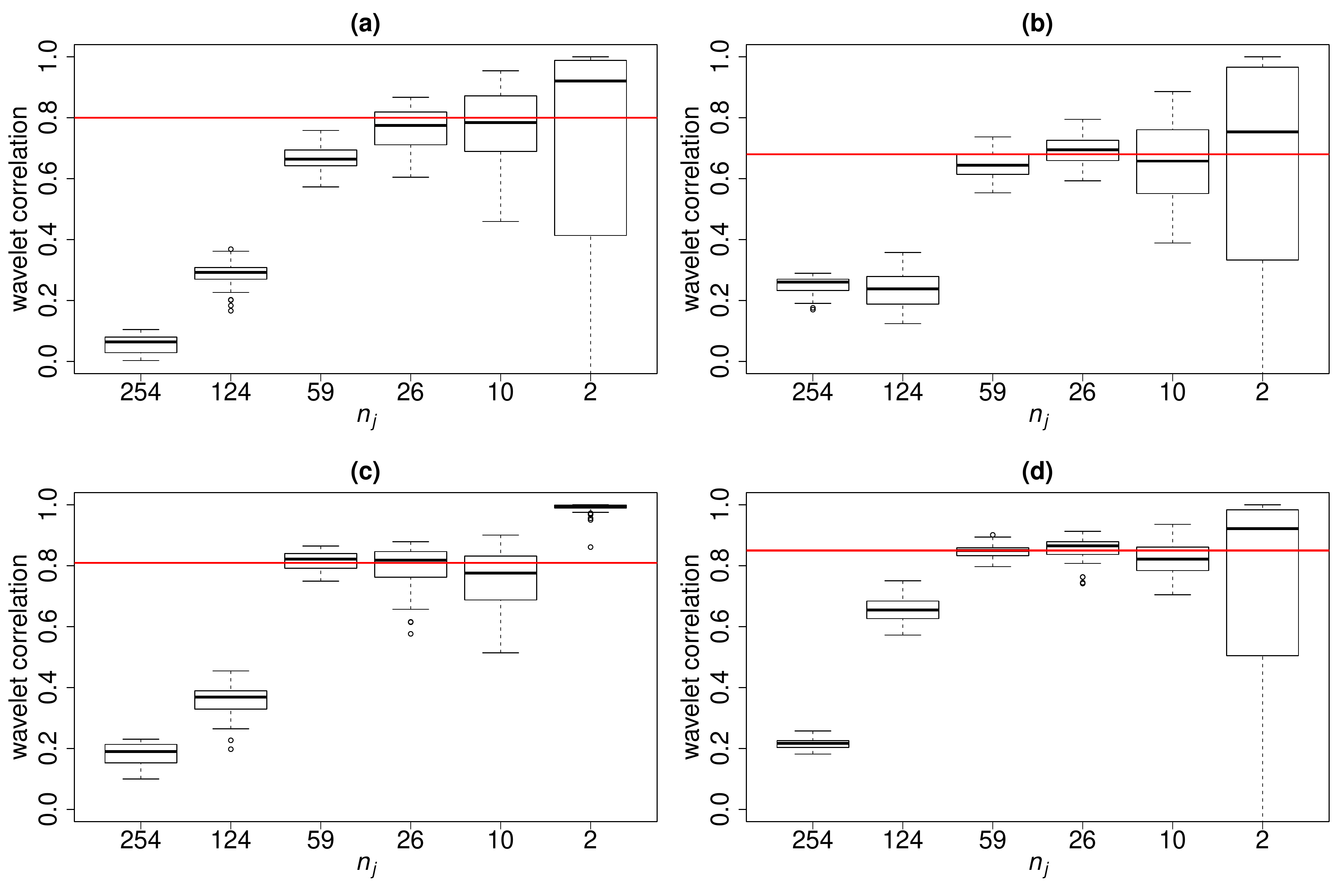}}
\end{center}
\end{figure}

As already shown in Figure \ref{fig:cov} for simulated data, the advantage of representing the wavelet correlation in terms of scale is to qualitatively determine the scales necessary to estimate the long-memory parameters and long-range covariance matrix. When dealing with real data, bootstrap is providing a way to assess the variability of the estimators. Using the real data described in Section \ref{sec:real_data_desc}, sliding overlapping window of the time series were extracted containing 512 points and we repeated the estimation until reaching the final point of the time series. This is illustrated in Figure \ref{fig:cov_real_data}, where an example of four pairs of fMRI data from one subject is presented. Boxplots are constructed using the sliding window extractions. From these plots, and taking into account %the length of the time series and 
neuroscientific hypothesis stating that the signal of interest for resting state is occurring for frequency below 0.1Hz, we chose to compute the long-memory parameters between scales 3 and 6. 

\begin{figure}[!ht]
\caption{Histogram of $\hat\bd$ from a subject of fMRI data set.}
\label{fig:hist}
\begin{center}
{\includegraphics[height=6cm]{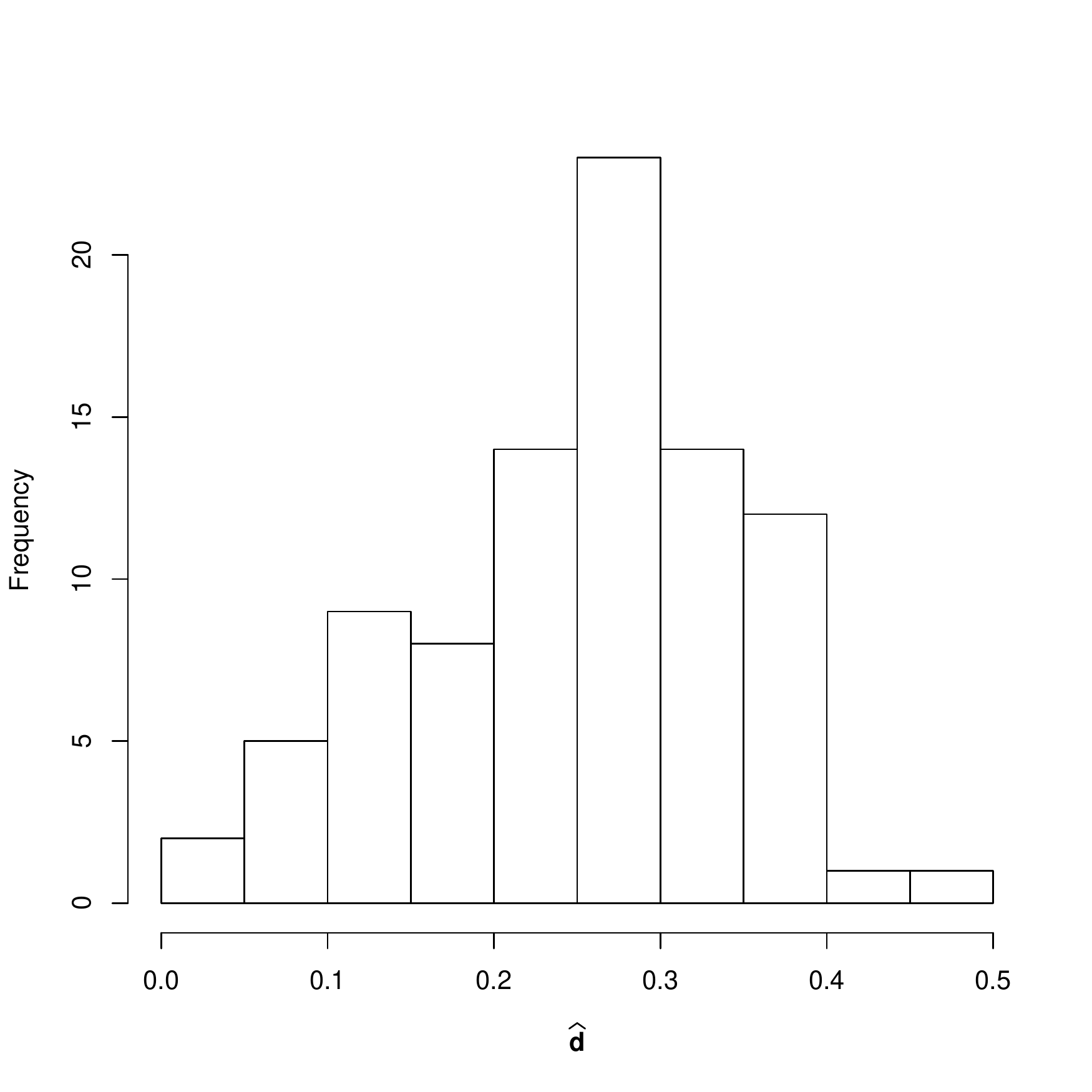}}
\end{center}
\end{figure}

\begin{figure}[!h]
\caption{Estimation of $\bOmega$ from a subject of fMRI data set.}
\label{fig:corr}
\begin{center}
{\includegraphics[height=9cm]{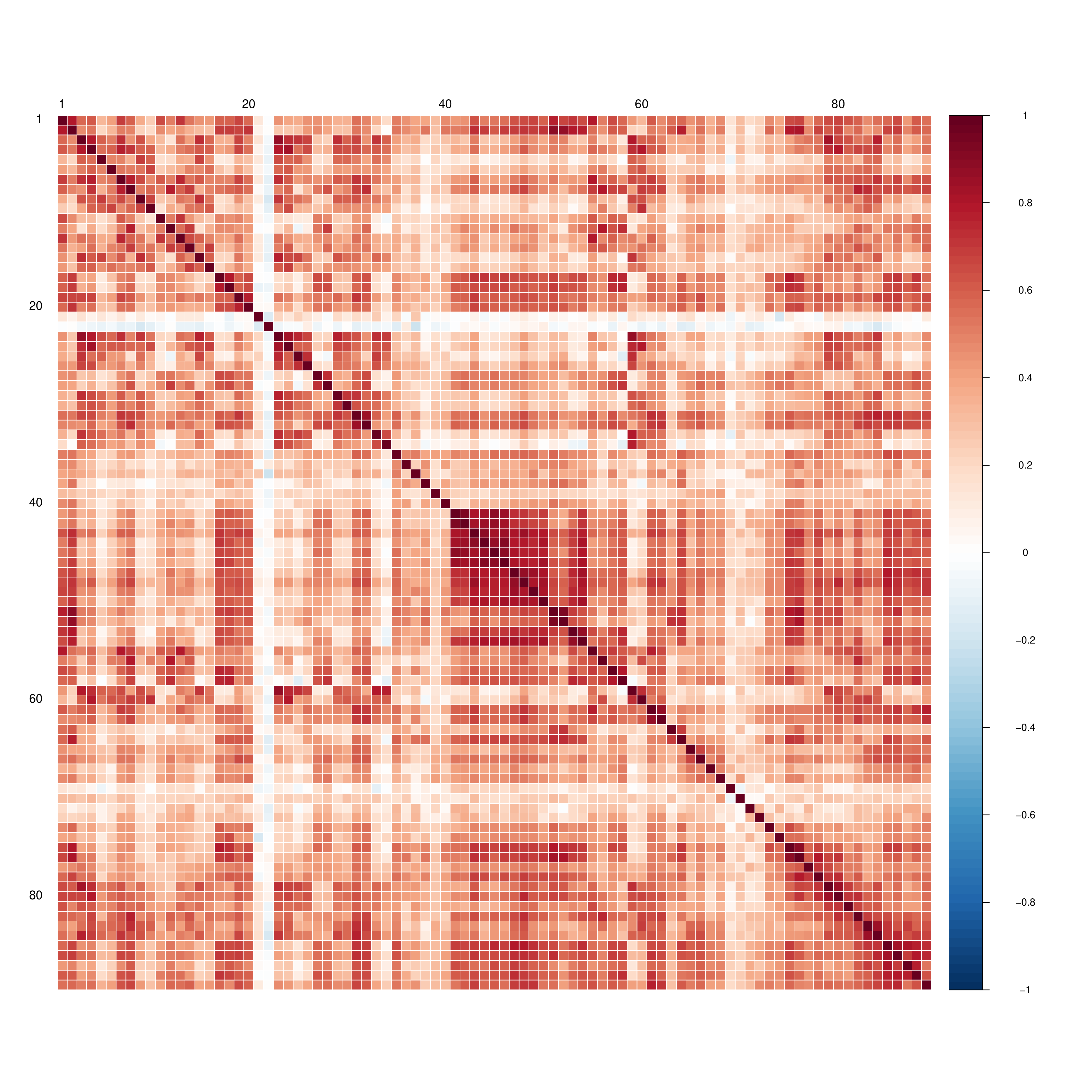}}
\end{center}
\end{figure}

Figures \ref{fig:hist} and \ref{fig:corr} display an example of long-memory parameter and long-run correlation estimated for one subject.


\newpage

\section*{Conclusion}

The \textsf{R} package \pkg{multiwave} provides a versatile wavelet-based approach, as well as a Fourier-based approach, for estimating long-memory parameters and long-run covariance matrices of multivariate time series. The two estimation procedures are based on semi-parametric approaches of \cite{Shimotsu07} and \cite{AchardGannaz2014}. The added value of the package is to provide estimations in long-range dependence multidimensional settings, which is not proposed presently by any \textsf{R} package to our knowledge. This paper describes the functions of the package \pkg{multiwave} and discusses some practical points for applications, including on a real data set. A simulation study shows first that multivariate estimation improves univariate estimation. The advantage of the wavelet-based procedure with respect to the Fourier-based estimation is its flexibility, allowing to take into account trends or nonstationarity.

\vspace{10pt}

{\bf Acknowledgements.} The authors are grateful to Shimotsu (\url{http://shimotsu.web.fc2.com/Site/Matlab\_Codes.html}) and to Fa\"y, Moulines, Roueff and Taqqu for kindly providing the codes of their respective papers. This work was partly supported by the project \textit{Graphsip} from Agence Nationale de la Recherche (ANR-14-CE27-0001). S.A. was partly funded by a grant from la Région Rhône-Alpes and a grant from
AGIR-PEPS, Université Grenoble Alpes–CNRS.

\bibliographystyle{plain}
\bibliography{biblio}

\newpage

\begin{appendices}

\pagestyle{empty}

\begin{table}
\caption{Multivariate wavelet Whittle estimation of $\bd$ for a bivariate $ARFIMA(0,\bd,0)$ with $\rho=0.8$, $N=512$ with 1000 repetitions. For the estimation, $j_0=1$. }
\label{tab:dwav}
\begin{center}
\begin{tabular}{@{} S@{~~}S@{~~}S@{~~}S@{~~}S@{~~}} 
\toprule
 $\bd$ & \multicolumn{1}{c}{\emph{bias}} & \multicolumn{1}{c}{\emph{std}} & \multicolumn{1}{c}{\emph{RMSE}} & \multicolumn{1}{c}{\emph{ratio M/U}} \\
\midrule 
 0.2  & -0.017 & 0.039 & 0.0425 & 0.7785 \\ 
 0.0   & 0.013 & 0.0391 & 0.0412 & 0.9014 \\  \hline
 0.2  & -0.0313 & 0.0376 & 0.049 & 0.8960 \\ 
 0.2  & -0.0316 & 0.0378 & 0.0493 & 0.8805 \\  \hline
 0.2  & -0.017 & 0.0383 & 0.0419 & 0.7673 \\ 
 0.4   & -0.0442 & 0.0395 & 0.0592 & 0.7902 \\
\bottomrule

\end{tabular}
\end{center}
\end{table}

\begin{table}
\caption{Wavelet Whittle estimation of $\bOmega$ for a bivariate $ARFIMA(0,\bd,0)$ with $\rho=0.8$, $N=512$ with 1000 repetitions. For the estimation, $j_0=1$.}
\label{tab:Owav}
\begin{center}
\begin{tabular}{@{}l@{}S@{}S@{}S@{}S@{~~~}S@{}S@{}S@{}S@{~~~}S@{}S@{}S@{}S} 
  \toprule
  & \multicolumn{3}{c}{\textbf{$\bd=( 0.2 , 0 )$ }} & & \multicolumn{3}{c}{\textbf{$\bd=( 0.2 , 0.2 )$ } }& & \multicolumn{3}{c}{\textbf{$\bd=( 0.2 , 0.4 )$ } }\\
  \cmidrule{2-4} \cmidrule{6-8} \cmidrule{10-12}
  
  $\bOmega$ &  \multicolumn{1}{r}{\emph{bias}} &  \multicolumn{1}{r}{\emph{std}} &  \multicolumn{1}{r}{\emph{RMSE}} & &   \multicolumn{1}{r}{\emph{bias}} &  \multicolumn{1}{r}{\emph{std}} &  \multicolumn{1}{r}{\emph{RMSE}} & &   \multicolumn{1}{r}{\emph{bias}} &  \multicolumn{1}{r}{\emph{std}} &  \multicolumn{1}{r}{\emph{RMSE}}\\
    \cmidrule{2-4} \cmidrule{6-8} \cmidrule{10-12}

  $\Omega_{1,1}$  & 0.0417 & 0.0724 & 0.0836 & & 0.0343 & 0.0711 & 0.0789 & & 0.0417 & 0.0717 & 0.083 \\ 
  $\Omega_{1,2}$  & 0.0382 & 0.0657 & 0.0759 & & 0.0279 & 0.0626 & 0.0686 & & 0.0673 & 0.0684 & 0.0959\\ 
  $\Omega_{2,2}$  & 0.0048 & 0.0709 & 0.071 & & 0.0323 & 0.0714 & 0.0784 && 0.0748 & 0.0748 & 0.1057\\ 
  correlation    &  0.0191 & 0.0227 & 0.0296 &    &  0.001 & 0.0164 & 0.0164 &   &  0.0194 & 0.0234 & 0.0304\\ 
 \bottomrule
\end{tabular}
\end{center}
\end{table}
\begin{table}
\caption{Multivariate wavelet Whittle estimation of $\bd$ for a bivariate $ARFIMA(1,\bd,0)$ with $\rho=0.8$, $N=512$ with 1000 repetitions. For the estimation, $j_0=3$.}
\label{tab:dwav2}
\begin{center}
\begin{tabular}{@{}S@{~~~}S@{~~~}S@{~~~}S@{~~~}S@{}} 
\toprule
  $\bd$ & \multicolumn{1}{c}{\emph{bias}} & \multicolumn{1}{c}{\emph{std}} & \multicolumn{1}{c}{\emph{RMSE}} & \multicolumn{1}{c}{\emph{ratio M/U}} \\
\midrule 
  0.2  & -0.0473 & 0.1213 & 0.1302 & 0.8472 \\ 
 0.0  & -0.0371 & 0.1266 & 0.132 & 0.8511 \\  \hline
  0.2  & -0.0623 & 0.1209 & 0.136 & 0.8848 \\ 
 0.2  & -0.0526 & 0.1258 & 0.1364 & 0.8714 \\  \hline
 0.2  & -0.066 & 0.1244 & 0.1408 & 0.9161 \\ 
 0.4   & -0.0584 & 0.1293 & 0.1418 & 0.8935 \\  
\bottomrule

\end{tabular}
\end{center}
\end{table}

\begin{table}[!ht]
  \caption{Wavelet Whittle estimation of $\bOmega$ for a bivariate $ARFIMA(1,\bd,0)$ with $\rho=0.8$, $N=512$ with 1000 repetitions. For the estimation, $j_0=3$.}
  \begin{center}
\label{tab:Owav2}
\begin{tabular}{@{}l@{}S@{}S@{}S@{}S@{~~~}S@{}S@{}S@{}S@{~~~}S@{}S@{}S@{}S} %\hspace{0.5cm}
  \toprule
  & \multicolumn{3}{c}{\textbf{$\bd=( 0.2 , 0 )$ }} & & \multicolumn{3}{c}{\textbf{$\bd=( 0.2 , 0.2 )$ } }& & \multicolumn{3}{c}{\textbf{$\bd=( 0.2 , 0.4 )$ } }\\
  \cmidrule{2-4} \cmidrule{6-8} \cmidrule{10-12}
  
  $\bOmega$ &  \multicolumn{1}{r}{\emph{bias}} &  \multicolumn{1}{r}{\emph{std}} &  \multicolumn{1}{r}{\emph{RMSE}} & &   \multicolumn{1}{r}{\emph{bias}} &  \multicolumn{1}{r}{\emph{std}} &  \multicolumn{1}{r}{\emph{RMSE}} & &   \multicolumn{1}{r}{\emph{bias}} &  \multicolumn{1}{r}{\emph{std}} &  \multicolumn{1}{r}{\emph{RMSE}}\\
    \cmidrule{2-4} \cmidrule{6-8} \cmidrule{10-12}

$\Omega_{1,1}$  & -0.0067 & 0.0828 & 0.0831 & & 0.0024 & 0.0859 & 0.0860 & & 0.005 & 0.0889 & 0.089 \\ 
  $\Omega_{1,2}$  & 0.0551 & 0.0828 & 0.0995 & & 0.0495 & 0.0794 & 0.0936 & & 0.0507 & 0.0875 & 0.1012 \\ 
  $\Omega_{2,2}$  & 0.1363 & 0.131 & 0.1891 & & 0.1412 & 0.1384 & 0.1977 & & 0.1391 & 0.1426 & 0.1992 \\ 
   correlation       &  -0.0088 & 0.0805 & 0.081 & & -0.0386 & 0.0527 & 0.0653 & &  -0.0358 & 0.0985 & 0.1047 \\ 
\bottomrule
\end{tabular}
\end{center}
\end{table}

\begin{table}
\caption{Multivariate wavelet Whittle estimation of $\bd$ for a bivariate $ARFIMA(0,\bd,0)$ with $\rho=0.8$, $N=512$ with 1000 repetitions. Nonstationary cases. For the estimation, $j_0=2$.}
\label{tab:dwav3}
\begin{center}
\begin{tabular}{c@{~~~~}c}
\begin{tabular}{@{} S@{~}S@{~}S@{~}S@{~}S@{}}
\toprule
 $\bd$ & \multicolumn{1}{c}{\emph{bias}} & \multicolumn{1}{c}{\emph{std}} & \multicolumn{1}{c}{\emph{RMSE}} & \multicolumn{1}{c}{\emph{ratio M/U}} \\
\midrule
1.2  & -0.0338 & 0.0762 & 0.0834 & 0.8510 \\ 
 1   & -0.0276 & 0.0725 & 0.0776 & 0.8316 \\  \hline
1.2  & -0.043 & 0.0732 & 0.0849 & 0.8672 \\ 
 1.2   & -0.0411 & 0.0743 & 0.0849 & 0.8591 \\  \hline
1.2  & -0.0338 & 0.0741 & 0.0814 & 0.8310 \\ 
 1.4  & -0.0356 & 0.0797 & 0.0873 & 0.8344 \\  
\bottomrule
\end{tabular}
&
\begin{tabular}{@{} S@{~}S@{~}S@{~}S@{~}S@{}} %\hspace{0.5cm}
\toprule
 $\bd$ & \multicolumn{1}{c}{\emph{bias}} & \multicolumn{1}{c}{\emph{std}} & \multicolumn{1}{c}{\emph{RMSE}} & \multicolumn{1}{c}{\emph{ratio M/U}} \\%& ratio M/U
\midrule
2.2  & -0.0421 & 0.0884 & 0.0979 & 0.8718 \\ 
 2  & -0.0403 & 0.0862 & 0.0951 & 0.8516 \\  \hline
2.2  & -0.0503 & 0.086 & 0.0996 & 0.8874 \\ 
 2.2  & -0.049 & 0.0823 & 0.0958 & 0.8566 \\  \hline
 2.2  & -0.0436 & 0.0868 & 0.0971 & 0.8651 \\ 
 2.4  & -0.0429 & 0.0831 & 0.0935 & 0.8400 \\ 
 
 \bottomrule

\end{tabular}
\end{tabular}
\end{center}
\end{table}

\begin{table}
\caption{Wavelet Whittle estimation of $\bOmega$ for a bivariate $ARFIMA(0,\bd,0)$ with $\rho=0.8$, $N=512$ with 1000 repetitions. Nonstationary cases. For the estimation, $j_0=2$.}
\label{tab:Owav3}
\begin{center}
\begin{tabular}{@{}l@{}S@{}S@{}S@{}S@{~~~}S@{}S@{}S@{}S@{~~~}S@{}S@{}S@{}S} %\hspace{0.5cm}
  \toprule
  & \multicolumn{3}{c}{\textbf{$\bd=( 1.2 , 1 )$ }} & & \multicolumn{3}{c}{\textbf{$\bd=( 1.2 , 1.2 )$ } }& & \multicolumn{3}{c}{\textbf{$\bd=( 1.2 , 1.4 )$ } }\\
  \cmidrule{2-4} \cmidrule{6-8} \cmidrule{10-12}
  
  $\bOmega$ &  \multicolumn{1}{r}{\emph{bias}} &  \multicolumn{1}{r}{\emph{std}} &  \multicolumn{1}{r}{\emph{RMSE}} & &   \multicolumn{1}{r}{\emph{bias}} &  \multicolumn{1}{r}{\emph{std}} &  \multicolumn{1}{r}{\emph{RMSE}} & &   \multicolumn{1}{r}{\emph{bias}} &  \multicolumn{1}{r}{\emph{std}} &  \multicolumn{1}{r}{\emph{RMSE}}\\
    \cmidrule{2-4} \cmidrule{6-8} \cmidrule{10-12}

$\Omega_{1,1}$  & -0.0049 & 0.1362 & 0.1363 & & 0.0059 & 0.1369 & 0.137 &&  -0.0047 & 0.1361 & 0.1361 \\ 
  $\Omega_{1,2}$  &  0.018 & 0.1168 & 0.1182 & & 0.0093 & 0.1155 & 0.1158 & & 0.0159 & 0.1266 & 0.1276 \\ 
  $\Omega_{2,2}$  & -0.0027 & 0.1277 & 0.1277 & & 0.0042 & 0.1386 & 0.1386  & & -0.0113 & 0.1487 & 0.1491\\ 
   correlation       &  0.0214 & 0.0475 & 0.0521  &   &  0.0051 & 0.0286 & 0.0291  &   &  0.0227 & 0.0507 & 0.0555\\  

\bottomrule
\end{tabular}

\vspace*{1cm}

\begin{tabular}{@{}l@{}S@{}S@{}S@{}S@{~~~}S@{}S@{}S@{}S@{~~~}S@{}S@{}S@{}S} %\hspace{0.5cm}
  \toprule
  & \multicolumn{3}{c}{\textbf{$\bd=( 2.2 , 2 )$ }} & & \multicolumn{3}{c}{\textbf{$\bd=( 2.2 , 2.2 )$ } }& & \multicolumn{3}{c}{\textbf{$\bd=( 2.2 , 2.4 )$ } }\\
  \cmidrule{2-4} \cmidrule{6-8} \cmidrule{10-12}
  
  $\bOmega$ &  \multicolumn{1}{r}{\emph{bias}} &  \multicolumn{1}{r}{\emph{std}} &  \multicolumn{1}{r}{\emph{RMSE}} & &   \multicolumn{1}{r}{\emph{bias}} &  \multicolumn{1}{r}{\emph{std}} &  \multicolumn{1}{r}{\emph{RMSE}} & &   \multicolumn{1}{r}{\emph{bias}} &  \multicolumn{1}{r}{\emph{std}} &  \multicolumn{1}{r}{\emph{RMSE}}\\
    \cmidrule{2-4} \cmidrule{6-8} \cmidrule{10-12}

$\Omega_{1,1}$  &  -0.0383 & 0.1795 & 0.1835 & & -0.0253 & 0.1789 & 0.1807 & & -0.0361 & 0.1776 & 0.1812\\ 
  $\Omega_{1,2}$  &  -0.0043 & 0.1565 & 0.1565 & & -0.0129 & 0.1493 & 0.1498& & -0.0097 & 0.1602 & 0.1605\\
  $\Omega_{2,2}$  &  -0.0318 & 0.1776 & 0.1804 & & -0.0276 & 0.1809 & 0.183 & & -0.0481 & 0.1813 & 0.1876 \\
   correlation       &  0.0251 & 0.0604 & 0.0654   &   &  0.0087 & 0.0374 & 0.0384  &   &  0.0249 & 0.0626 & 0.0674\\ 

\bottomrule
\end{tabular}

\end{center}
\end{table}

\begin{table}
\caption{Multivariate Fourier Whittle estimation of $\bd$ for a bivariate $ARFIMA(0,\bd,0)$ with $\rho=0.8$, $N=512$ with 1000 repetitions depending on the number of frequencies $m$. $\lfloor x\rfloor$ denotes the closest integer smaller than $x$.}
\label{tab:dfou}
\begin{center}
\begin{tabular}{r@{~} S@{~}S@{~}S@{~}S@{~}S@{~~}S@{~}S@{~}S@{~}S@{~}S@{}} %\hspace{0.5cm}
\toprule
& \multicolumn{4}{c}{\textbf{$m=\lfloor N^{0.65}\rfloor =57$}} & & \multicolumn{5}{c}{\textbf{$m=\lfloor N^{\eta}\rfloor$} }\\
\cmidrule{2-5} \cmidrule{7-11}
$\bd$ & \multicolumn{1}{c}{\emph{bias}} & \multicolumn{1}{r}{\emph{std}} & \multicolumn{1}{r}{\emph{RMSE}} & \multicolumn{1}{r}{\emph{W/F}} & & $\eta$& \multicolumn{1}{r}{\emph{bias}} & \multicolumn{1}{r}{\emph{std}} & \multicolumn{1}{r}{\emph{RMSE}} & \multicolumn{1}{r}{\emph{W/F}}\\
\midrule 
  0.2 & -0.002 & 0.0576 & 0.047 & 0.9049 & & 0.90 & -0.0197 & 0.0271 & 0.0335 & 1.2689\\ 
 0  & -9e-04 & 0.0593 & 0.0238 & 1.7295  & &    & -0.0033 & 0.0264 & 0.0267 & 1.5475\\ \hline
 0.2 & -0.0033 & 0.0574 & 0.0531 & 0.9214 & & 0.85 & -0.013 & 0.0306 & 0.0332 & 1.4750\\ 
 0.2  & -0.0031 & 0.0591 & 0.0522 & 0.9434  &&       & -0.0123 & 0.0293 & 0.0318 & 1.5503\\ \hline
 0.2 & 8e-04 & 0.0576 & 0.0579 & 0.7239 && 0.85  & -0.0136 & 0.0308 & 0.0337 & 1.2455\\ 
 0.4  & 9e-04 & 0.0595 & 0.0889 & 0.6666 &&       &-0.0192 & 0.0299 & 0.0355 & 1.6690\\ 
\bottomrule
\end{tabular}
\end{center}

 \end{table}

\thispagestyle{empty}
\begin{table}
\caption{Fourier Whittle estimation of $\bOmega$ for a bivariate $ARFIMA(0,\bd,0)$ with $\rho=0.8$, $N=512$ with 1000 repetitions. The number of frequencies is $m=\lfloor n^{\eta}\rfloor$ with $\eta=0.65$ as chosen in \cite{Shimotsu07}.}
\label{tab:Ofou}
\begin{center}

\begin{tabular}{@{}l@{}S@{}S@{}S@{}S@{~}S@{}S@{}S@{}S@{~}S@{}S@{}S@{}S@{}} %\hspace{0.5cm}
  \toprule
  & \multicolumn{3}{c}{\textbf{$\bd=( 0.2 , 0 )$ }} & & \multicolumn{3}{c}{\textbf{$\bd=( 0.2 , 0.2 )$ } }& & \multicolumn{3}{c}{\textbf{$\bd=( 0.2 , 0.4 )$ } }\\
  \cmidrule{2-4} \cmidrule{6-8} \cmidrule{10-12}
  
  $\bOmega$ &  \multicolumn{1}{r}{\emph{bias}}  &  \multicolumn{1}{r}{\emph{RMSE}} & \multicolumn{1}{r}{\emph{W/F}} & &   \multicolumn{1}{c}{\emph{bias}} &  \multicolumn{1}{r}{\emph{RMSE}} & \multicolumn{1}{r}{\emph{W/F}} & &   \multicolumn{1}{r}{\emph{bias}}  &  \multicolumn{1}{r}{\emph{RMSE}} & \multicolumn{1}{r}{\emph{W/F}}\\
  \cmidrule{2-4} \cmidrule{6-8} \cmidrule{10-12}

$\Omega_{1,1}$  & 0.0186& 0.2027 & 0.4123 & & 0.022 & 0.2036 & 0.3876  &  & 0.0113  & 0.201 & 0.4127 \\ 
  $\Omega_{1,2}$ & 0.0124 & 0.0759 & 0.4381 & & 0.0197  & 0.0686 & 0.3923 & & 0.0149  & 0.0959 & 0.5545 \\ 
  $\Omega_{2,2}$  & 0.0146  & 0.2138 & 0.3322 & & 0.0263  & 0.2163 & 0.3624 & & 0.0291  & 0.2162 & 0.4891 \\ 
   correlation   &  -0.0014  & 0.0357 & 0.8298  &   &  -2e-04  & 0.0355 & 0.4621  &   &  -0.0018  & 0.0359 & 0.8483\\  

\bottomrule
\end{tabular}
\end{center}
\end{table}

\begin{table}
\caption{Fourier Whittle estimation of $\bOmega$ for a bivariate $ARFIMA(0,\bd,0)$ with $\rho=0.8$, $N=512$ with 1000 repetitions. The number of frequencies is $m=\lfloor n^{\eta}\rfloor$ with $\eta$ such that RMSE of $\hat\bd$ is minimized.}

\label{tab:Ofou2}
\begin{center}

\begin{tabular}{@{}l@{}S@{}S@{}S@{}S@{~}S@{}S@{}S@{}S@{~}S@{}S@{}S@{}S@{}} %\hspace{0.5cm}
  \toprule
  & \multicolumn{3}{c}{\textbf{$\bd=( 0.2 , 0 )$ }} & & \multicolumn{3}{c}{\textbf{$\bd=( 0.2 , 0.2 )$ } }& & \multicolumn{3}{c}{\textbf{$\bd=( 0.2 , 0.4 )$ } }\\
   & \multicolumn{3}{c}{$\eta=0.9$} & & \multicolumn{3}{c}{$\eta=0.85$} & & \multicolumn{3}{c}{$\eta=0.85$}\\
  \cmidrule{2-4} \cmidrule{6-8} \cmidrule{10-12}
  
  $\bOmega$ &  \multicolumn{1}{r}{\emph{bias}}  &  \multicolumn{1}{r}{\emph{RMSE}} & \multicolumn{1}{r}{\emph{W/F}} & &   \multicolumn{1}{r}{\emph{bias}} &  \multicolumn{1}{r}{\emph{RMSE}} & \multicolumn{1}{r}{\emph{W/F}} & &   \multicolumn{1}{r}{\emph{bias}}  &  \multicolumn{1}{r}{\emph{RMSE}} & \multicolumn{1}{r}{\emph{W/F}}\\
  \cmidrule{2-4} \cmidrule{6-8} \cmidrule{10-12}

$\Omega_{1,1}$  & 0.0622  & 0.094 & 0.8896 && 0.0387 & 0.0832 & 0.9484  && 0.0389  & 0.0833 & 0.9962  \\ 
  $\Omega_{1,2}$ & 0.0222  & 0.0638 & 1.1907 & & 0.0304 & 0.0731 & 0.9379 & & 0.047 & 0.0826 & 1.1615 \\ 
  $\Omega_{2,2}$  & -0.0031 & 0.0637 & 1.1149  & & 0.0373  & 0.0839 & 0.9341  & & 0.0812  & 0.1132 & 0.9338 \\ 
   correlation   &  -0.0013  & 0.0163 & 1.8222  &   &  -4e-04  & 0.0177 & 0.9237    &   &  -0.0012  & 0.0179 & 1.6969  \\  

\bottomrule
\end{tabular}
\end{center}
\end{table}

\end{appendices}

\end{document}